\documentclass{article}
\usepackage{hyperref}
\usepackage[american]{babel}
\usepackage{amsfonts,amsmath,amssymb,epsf,epsfig}

\usepackage{amsthm}

\usepackage{xypic}
\xyoption{all}
\input{xypic}
\newdir{ >}{{}*!/-8pt/\dir{>}}

\newcommand{\id}{{\rm id}} 
\renewcommand{\:}{\colon} 
 
\renewcommand{\hat}{\widehat} 

\newcommand{\Hom}{\mathop{{\rm Hom}}\nolimits}

\def\R{\mathbb{R}}

\def\C{\mathbb{C}}

\def\Q{\mathbb{Q}}

 \def\ad{{\rm ad}}     \def\Hom{{\rm Hom}}
\def\id{{\rm id}}

\newtheorem{theo}{Theorem}[section]
\newtheorem{Defi}[theo]{{\bf Definition}}
\newenvironment{defi}{\begin{Defi} \normalfont}{\end{Defi}}

\newtheorem{Prop}[theo]{{\bf Proposition}}
\newenvironment{prop}{\begin{Prop} \normalfont}{\end{Prop}}
\newtheorem{Cor}[theo]{{\bf Corollary}}
\newenvironment{cor}{\begin{Cor} \normalfont}{\end{Cor}}
\newtheorem{Lem}[theo]{{\bf Lemma}}
\newenvironment{lem}{\begin{Lem} \normalfont}{\end{Lem}}
\newtheorem{pro}[theo]{Problem}
\newtheorem{Exa}[theo]{{\bf Example}}

\newtheorem{Exem}{{\bf Example}}[section]
\newenvironment{exem}{\begin{Exem} \strut\\ 
\normalfont}{\hfill$\Diamond$\end{Exem}}

\newtheorem{Rem}[theo]{{\bf Remark}}
\newenvironment{rem}{\begin{Rem} \normalfont}{\end{Rem}}

\newenvironment{prf}{\begin{proof}}{\end{proof}}

\begin{document}
\title{Algebraic deformation quantization of Leibniz algebras} 

\author{Charles Alexandre \\
Universit\'{e} de Strasbourg
  \thanks{Char.Alexandre@gmail.com}\\
\and Martin Bordemann \\
Universit\'{e} de Mulhouse
  \thanks{Martin.Bordemann@uha.fr}\\
\and Salim Rivi\`ere\\
Universit\'{e} d'Angers
\thanks{salim.riviere@univ-angers.fr} \\
\and Friedrich Wagemann\\
     Universit\'{e} de Nantes
     \thanks{wagemann@math.univ-nantes.fr}}   
\maketitle

\begin{abstract}
 In this paper we focus on a certain self-distributive multiplication on
 coalgebras, which leads to so-called rack bialgebra. We construct canonical 
 rack bialgebras (some kind of enveloping algebras) for any Leibniz algebra.  
 
Our motivation is deformation quantization of Leibniz algebras in the sense 
of \cite{DW2013}.  
 Namely, the canonical rack bialgebras we have constructed for any Leibniz 
algebra
 lead to a simple explicit formula of the rack-star-product
 on the dual of a Leibniz algebra recently constructed by Dherin and 
 Wagemann in \cite{DW2013}. We clarify this framework setting up a general 
deformation theory for rack bialgebras and show that the rack-star-product
turns out to be a deformation of the trivial rack bialgebra product.  
\end{abstract}

\section*{Introduction}

\subsection*{The algebraic structures involved in Leibniz deformation 
quantization}

Recall that a \emph{pointed rack} (see \cite{FR92})
is a pointed set $(X,e)$ together with a binary operation 
$\rhd:X\times X\to X$ such that for all $x\in X$, the map $y\mapsto x\rhd y$ 
is bijective and such that for all $x,y,z\in X$,
the self-distributivity and unit relations
$$
   x\rhd(y\rhd z)\,=\,(x\rhd y)\rhd(x\rhd z),~~~e\triangleright x = x,~~~
   \mathrm{and}~~~x\triangleright e = e
$$
are satisfied. Observe that racks are not algebras over an operad, but the 
correct algebraic structure is that of a properad. Therefore the standard 
deformation 
theory of algebras over an operad does not apply.  
Imitating the notion of a Lie group, the smooth version of a pointed rack 
is called a {\it Lie rack}. 

An important class of examples of racks are the so-called 
{\it augmented racks}, 
see \cite{FR92}. 
An augmented rack is the data of a group $G$, a $G$-set $X$ and a 
map $p:X\to G$ 
such that for all $x\in X$ and all $g\in G$, 
$$p(g\cdot x)\,=\,gp(x)g^{-1}.$$
The set $X$ becomes then a rack by setting $x\rhd y\,:=\,p(x)\cdot y$. 
In fact, augmented racks are the Drinfeld center (or the Yetter-Drinfeld 
modules)
in the monoidal category of $G$-sets over the (set-theoretical) Hopf 
algebra $G$, see for example
\cite{KW14}. Any rack may be augmented in many ways, for example by
 using the canonical 
morphism to its associated group (see \cite{FR92}) or to its group of 
bijections or 
to its group of automorphisms.   

In order to 
formalize the notion of a rack, one needs the diagonal map 
${\rm diag}_X:X\to X\times X$ given by $x\mapsto(x,x)$. Then the
self-distributivity relation reads in terms of maps  
\begin{eqnarray*}
\lefteqn{\mathbf{m}\circ (\mathrm{id}_M \times \mathbf{m})}  \\
  & =  & \mathbf{m}\circ (\mathbf{m}\times \mathbf{m})\circ 
     (\mathrm{id}_M\times\tau_{M,M}\times \mathrm{id}_M) \circ 
     (\mathrm{diag}_M\times\mathrm{id}_M\times \mathrm{id}_M). 
\end{eqnarray*}
Axiomatizing this kind of structure, one may start with a coalgebra $C$ 
and look for rack operations on this fixed coalgebra, see \cite{CCES08}
and \cite{L12}.
A natural framework where this kind of structure arises 
is by taking point-distributions 
over (resp. to) the pointed manifold given by a Lie rack, see \cite{Ser64},
\cite{Bou}, \cite{Mos} or \cite{ABRW}. 
We dub the arising structure a {\it rack bialgebra}, see Definition
\ref{definition_rack_bialgebra}. We carry out some structure theory for 
rack bialgebras based on semigroup theory in the article \cite{ABRW}.  
  
Lie racks are intimately related to
{\it Leibniz algebras} ${\mathfrak h}$, i.e. a vector space ${\mathfrak h}$ 
with a 
bilinear bracket $[,]:{\mathfrak h}\otimes{\mathfrak h}\to{\mathfrak h}$ 
such that
for all $X,Y,Z\in{\mathfrak h}$, $[X,-]$ acts as a derivation:
\begin{equation}    \label{LeibnizIdentity}
[X,[Y,Z]]\,=\,[[X,Y],Z]+[Y,[X,Z]].
\end{equation}
Indeed, Kinyon showed in \cite{Kin07} that the tangent space at $e\in H$ 
of a Lie rack $H$ carries 
a natural structure of a Leibniz algebra, generalizing the relation between 
a Lie group 
and its tangent Lie algebra. Conversely, every (finite dimensional real or 
complex) 
Leibniz algebra ${\mathfrak h}$
may be integrated into a Lie rack $R_{\mathfrak h}$ (with underlying manifold 
${\mathfrak h}$) 
using the rack product
\begin{equation}   \label{rack_structure}
X\blacktriangleright Y\,:=\,e^{{\rm ad}_X}(Y),
\end{equation}
noting that the exponential of the inner derivation ${\rm ad}_X$ for each 
$X\in{\mathfrak h}$
is an automorphism.   

\subsection*{Leibniz deformation quantization}
 
Given a finite-dimensional real Lie algebra $({\mathfrak g},[,])$, 
its dual vector space ${\mathfrak g}^*$ is a smooth manifold
which carries a Poisson bracket on its space of smooth functions, defined 
for all 
$f,g\in{\mathcal C}^{\infty}({\mathfrak g}^*)$ and all $\xi\in{\mathfrak g}^*$
by the {\it Kostant-Kirillov-Souriau formula}
$$\{f,g\}(\xi)\,:=\,\langle \xi,[df(\xi),dg(\xi)]\rangle.$$
Here $df(\xi)$ and $dg(\xi)$ are linear functionals on ${\mathfrak g}^*$, 
identified with elements of
${\mathfrak g}$.

In the same way, a general finite dimensional Leibniz algebra ${\mathfrak h}$ 
gives rise to a smooth manifold 
${\mathfrak h}^*$, which carries now some kind of generalized Poisson bracket, 
namely
$$\{f,g\}(\xi)\,:=\,-\langle \xi,[df(0),dg(\xi)]\rangle,$$
see \cite{DW2013} for an explanation why we believe that this is the 
correct bracket
to consider. In particular, this generalized Poisson bracket need not be 
skew-symmetric.
 
The quantization procedure of this generalized Poisson bracket 
proposed in \cite{DW2013} works as follows: 
The cotangent lift of the above rack product
$$X\blacktriangleright Y\,=\,e^{{\rm ad}_X}(Y)$$
is interpreted as a symplectic micromorphism. 
The generating function of this micromorphism serves then as a phase function 
in a Fourier integral operator, whose asymptotic expansion gives rise to a 
star-product. 

One main goal of the present article is to set up a purely algebraic framework 
in which one may deformation quantize the dual of a Leibniz algebra.
The main feature will be to recover --in a rather explicit algebraic manner--
the star-product which has been constructed in \cite{DW2013} by analytic 
methods, see Corollary \ref{main_result}. The explicit formula reads:

Let $f,g\in{\mathcal C}^{\infty}({\mathfrak h}^*)$.
\begin{equation}  
(f\rhd_{\hbar}g)(\alpha)\,=\,\sum_{r=0}^{\infty}\frac{\hbar^r}{r!}
\sum_{i_1,\ldots,i_r=1}^n\frac{\partial^r f}{\partial\alpha_{i_1}\dots
\partial\alpha_{i_r}}(0)\Big((\widetilde{\rm ad}_{i_1}\circ\ldots\circ
\widetilde{\rm ad}_{i_r})(g)\Big)(\alpha),
\end{equation} 
where we have chosen a basis in the finite dimensional Leibniz algebra 
${\mathfrak h}$ and the first order differential operators ${\rm ad}_{i_1}$  
are defined by
$$
(\widetilde{\rm ad}_i(f))(\alpha)\,=\,\sum_{j,k=1}^n\alpha_k\,c_{ij}^k\,
\frac{\partial f}{\partial\alpha_j}(\alpha),
$$  
where $c_{ij}^k$ are the structure constants of ${\mathfrak h}$ w.r.t. this 
basis. Our first main result is thus Corollary \ref{main_result} where we 
show that the above $\rhd_{\hbar}$ is indeed a rack product and that it is equal 
(up to a sign) to the one constructed in \cite{DW2013} by analytic 
methods. 

In the remaining part of the paper, we answer the natural question in which 
cohomology theory the first term of the above deformation quantization 
appears as a 2-cocycle. 
For this, we set up a general deformation theory 
framework in which the above star-product appears as a formal deformation, 
its infinitesimal term defining a 2-cocycle and thus 
a second cohomology class. The main result is here Theorem 
\ref{rack_bialgebra_cohomology_complex} where we show that a natural 
differential $d$ on the adjoint rack bialgebra complex satisfies indeed
$d^2=0$. This is a combinatorially involved computation. Observe that this 
cohomology theory is thus well-defined in all degrees, in contrast to the 
related cohomology theory in \cite{CCES08} for which only degrees up to 
$3$ exist for the moment. One main point in this part of the paper is that 
our deformation complex for rack bialgebras replaces the (non existing) 
rack cohomology complex with adjoint coefficients (cf the case of group 
cohomology for a group $G$ where cohomology with adjoint coefficients in 
$G$ does not exist, while the Hochschild cohomology of the group algebra 
$K[G]$ with values in $K[G]$ may play this role).
\\

{\bf Acknowledgements:}
F.W. thanks Universit\'e de Haute Alsace for an invitation during which the 
shape of this research project was defined. 
Some part of the results of this paper constitute the master thesis of C.A..
We all thank the referee for useful remarks leading to an improvement of the
manuscript. 

\section{Preliminaries}

\subsection{Rack bialgebras}

In the following, let $K$ be an associative commutative unital ring
containing all the rational numbers. In the main part of this paper,
we will assume $K=\R$ or $\C$. The symbol
$\otimes$ will always denote the tensor product of $K$-modules over
$K$. For any coalgebra $(C,\Delta)$ over $K$, we shall use 
Sweedler's notation $\Delta(a)=\sum_{(a)}a^{(1)}\otimes a^{(2)}$
for any $a\in A$. We will feel free to suppress the sum-sign in Sweedler's notation
in complicated formulas for typographical reasons. 

The following sections will all deal with the following type of 
\emph{nonassociative bialgebra}: Let
$(B,\Delta,\epsilon,\mathbf{1},\mu)$ be a $K$-module such that
$(B,\Delta,\epsilon,\mathbf{1})$ is a \emph{coassociative counital coaugmented
coalgebra} (a $C^3$-coalgebra), and such that the linear map 
$\mu:B\otimes B\to B$
(the multiplication) is a morphism of $C^3$-coalgebras (it satisfies in 
particular 
$\mu(\mathbf{1}\otimes\mathbf{1})=\mathbf{1}$). We shall call this situation
a \emph{nonassociative $C^3I$-bialgebra} (where $I$ stands for $\mathbf{1}$
being an \emph{idempotent} for the multiplication $\mu$). For another
nonassociative $C^3I$-bialgebra $(B',\Delta',\epsilon',\mathbf{1}',\mu')$
a $K$-linear map $\phi:B\to B'$ will be called a \emph{morphism of
nonassociative $C^3I$-bialgebras} iff it is a morphism of $C^3$-coalgebras
and is multiplicative in the usual sense 
$\phi\big(\mu(a\otimes b)\big)=\mu'\big(\phi(a)\otimes\phi(b)\big)$
for all $a,b\in B$.

\begin{defi} \label{definition_rack_bialgebra}
  A \textbf{rack bialgebra} $(B,\Delta,\epsilon,\mathbf{1},\mu)$ 
  is a nonassociative
  $C^3I$-bialgebra (where we write
  for all $a,b\in B$
  $\mu(a\otimes b)=:a\triangleright b$) such that the following identities
  hold for all $a,b,c\in B$
  \begin{eqnarray}
   \mathbf{1}\triangleright a &=& a,  \label{EqRackBialgOnePreservesAll}\\
   a\triangleright \mathbf{1} & = & \epsilon(a)\mathbf{1}, \label{EqRackBialgAllPreserveOne} \\
   a \triangleright (b\triangleright c) &=&
     \sum_{(a)} (a^{(1)}\triangleright b)\triangleright
               (a^{(2)}\triangleright c). \label{EqDefAlgSelfDistributivity}
 \end{eqnarray}
 The last condition (\ref{EqDefAlgSelfDistributivity}) is called the 
 \textbf{self-distributivity} condition.
\end{defi}
Note that we do not demand that the $C^3$-coalgebra $B$ should be cocommutative
nor connected.  Similar definitions have been proposed in \cite{CCES08} and
in \cite{L12}.

\begin{exem}
Any $C^3$ coalgebra $(C,\Delta,\epsilon,\mathbf{1})$ carries
a \emph{trivial rack bialgebra structure} defined by the left-trivial 
multiplicaton
\begin{equation}
    a \triangleright_0 b := \epsilon(a)b
\end{equation}
which in addition is easily seen to be associative and left-unital, but in general
not unital.
\end{exem}

Another method of constructing rack bialgebras is \emph{gauging}:
Let $(B,\Delta,\epsilon,\mathbf{1},\mu)$ a rack bialgebra --where we write 
$\mu(a\otimes b)=a\triangleright b$ for all $a,b\in B$--, and let
$f:B\to B$ a morphism of $C^3$-coalgebras such that for all $a,b\in B$
\begin{equation}
   f(a\triangleright b)= a\triangleright \big(f(b)\big),
\end{equation}
i.e. $f$ is $\mu$-equivariant. 
It is a routine check that $(B,\Delta,\epsilon,\mathbf{1},\mu_f)$ is a 
rack bialgebra where for all $a,b\in B$ the multiplication is defined by
\begin{equation}
   \mu_f(a\otimes b):=a \triangleright_f b := \big(f(a)\big)\triangleright b.
\end{equation}
We shall call $(B,\Delta,\epsilon,\mathbf{1},\mu_f)$ the 
$f$-\emph{gauge
of $(B,\Delta,\epsilon,\mathbf{1},\mu)$}.

\begin{exem}
Let $(H,\Delta_H,\epsilon_H,\mu_H,\mathbf{1}_H,S)$ be a 
cocommutative Hopf algebra
over $K$. Then it is easy to see (cf.~also the particular case $B=H$ and 
$\Phi=\mathrm{id}_H$ of 
Proposition \ref{PAugmentedRackBialgebraIsRackBialgebra} for a detailed proof) 
that the new multiplication 
$\mu:H\otimes H\to H$, written $\mu(h\otimes h')=h\triangleright h'$,
defined by the usual \emph{adjoint representation}
\begin{equation}\label{EqSecondDefAdjointRep}
    h\triangleright h' := \mathrm{ad}_h(h')
       :=\sum_{(h)}h^{(1)}h'\big(S(h^{(2)})\big),
\end{equation}
equips the $C^4$-coalgebra $(H,\Delta_H,\epsilon_H,\mathbf{1}_H)$
with a rack bialgebra structure.
\end{exem}

In general, the adjoint representation does not seem to preserve the 
coalgebra structure if no cocommutativity is
assumed.

\begin{exem}  \label{rack_algebra}
Recall that a pointed set $(X,e)$ is a {\it pointed rack} in case there is a
binary operation 
$\rhd:X\times X\to X$ such that for all $x\in X$, the map $y\mapsto x\rhd y$ 
is bijective and such that for all $x,y,z\in X$,
the self-distributivity and unit relations
$$
   x\rhd(y\rhd z)\,=\,(x\rhd y)\rhd(x\rhd z),~~~e\triangleright x = x,~~~
   \mathrm{and}~~~x\triangleright e = e
$$
are satisfied. Then there is a natural rack bialgebra structure on the vector space
$K[X]$ which has the elements of $X$ as a basis. $K[X]$ carries the usual coalgebra 
structure such that all $x\in X$ are set-like: $\triangle(x)=x\otimes x$ for all 
$x\in X$. The product $\mu$ is then induced by the rack product. 
By functoriality, $\mu$ is compatible with $\triangle$ and $e$. 

Observe that this construction differs slightly from the construction in 
\cite{CCES08}, Section 3.1. 
\end{exem} 

\begin{rem}
It is shown in Theorem 4.3 of \cite{CCES08} that for a $C^4$-coalgebra $B$ with a 
self-distributive map $\lhd=q:B\otimes B\to B$ which is a morphism of coalgebras, the map 
$$
  R_q\,=\,({\rm id}_B\otimes q)\circ(\tau\otimes{\rm id}_B)\circ
  ({\rm id}_B\otimes\Delta)
$$
is a solution of the Yang-Baxter equation. We draw our reader's attention to the fact 
that Carter-Crans-Elhamdadi-Saito work in \cite{CCES08} with right racks, while we work here
with left racks. The statement of their theorem works also for left racks, but then one has 
to take
$$
   \tilde{R}_q\,=\,({\rm id}_B\otimes q)\circ(\Delta\otimes {\rm id}_B)\circ
    \tau.
$$   
In particular for any rack bialgebra,
$\tilde{R}_q$ is a solution of the Yang-Baxter equation.  
\end{rem}

\begin{exem} 
Here we suppose $K=\R$ or $\C$. 
Another general construction mechanism for rack bialgebras is exhibited in
\cite{ABRW}. Namely, let $(M,e,\rhd)$ be a Lie rack. In particular, $(M,e)$
is a pointed manifold and it makes sense to associate to it the vector space
$\mathcal{E}'_e(M)$ of distributions on $M$ which have their support in 
$\{e\}$. There is a corresponding functor 
$F:\mathcal{M}f*\to \mathbb{K}\mathbf{Vect}$ from the category of pointed 
manifolds $\mathcal{M}f*$ with values in the category of $K$-vector spaces 
$\mathbb{K}\mathbf{Vect}$. Observe that J.-P.
Serre \cite{Ser64} used it to define the universal enveloping algebra 
of a Lie group, see also \cite{Bou}, \cite{Mos} or \cite{ABRW}.

In fact, this framework is a special case of a local multiplication where
$0\times 0=0$. Such a multiplication gives rise to a bialgebra of
distributions supported at $0$. If the multiplication satisfies a certain
identity, the bialgebra satisfies the linearized identity. In our case,
the rack identity (i.e. self-distributivity) leads to the linearized
self-distributivity relation, i.e. condition
(\ref{EqDefAlgSelfDistributivity}). The general framework is described in
\cite{Mos} Section 3, see also \cite{PerShe}.  

In our case, the pointed manifold given by a 
Lie rack yields a rack bialgebra, see \cite{ABRW} where this functor is 
studied in detail.
\end{exem}    
 
We will also need the following structure:
\begin{defi}
 An \textbf{augmented rack bialgebra} over $K$ is a quadruple
 $(B,\Phi,H,\ell)$ consisting
 of a $C^3$-coalgebra $(B,\Delta,\epsilon,\mathbf{1})$, of a 
 cocommutative (!)
 Hopf algebra
 $(H,\Delta_H,\epsilon_H,\mathbf{1}_H,\mu_H,S)$, of a morphism of
 $C^3$-coalgebras $\Phi:B\to H$, and of a left action 
 $\ell:H \otimes B\to B$ of $H$ on $B$ which is a morphism of 
 $C^3$-coalgebras (i.e. $B$ is a $H$-module-coalgebra) 
 such that for all $h\in H$ and $a\in B$
 \begin{eqnarray}
     h.\mathbf{1} & = &\epsilon_H(h)\mathbf{1} 
            \label{EqDefAugBialgHPreservesE} \\
    \Phi(h.a) & = & \mathrm{ad}_h\big(\Phi(a)\big)
         \label{EqDefAugBialgPhiIntertwinesActionWithAdjoint}.
 \end{eqnarray}
 where $\mathrm{ad}$ denotes the usual adjoint representation for  
 Hopf algebras, see e.g. eqn (\ref{EqSecondDefAdjointRep}).
 
 \noindent We shall define a \textbf{morphism} 
 $(B,\Phi,H,\ell)\to (B',\Phi',H',\ell')$
\textbf{of augmented rack bialgebras} to be a pair $(\phi,\psi)$ 
of $K$-linear maps
where $\phi:(B,\Delta,\epsilon,\mathbf{1})\to (B',\Delta',\epsilon',\mathbf{1}')$
is a morphism of $C^3$-coalgebras, and $\psi:H\to H'$ is a morphism of Hopf 
algebras such that the obvious diagrams commute:
\begin{equation}\label{EqDefMorphAugmentedRackBialg}
   \Phi'\circ \phi = \psi\circ \Phi,~~~\mathrm{and}~~
         \ell'\circ (\psi\otimes \phi) = \phi\circ \ell 
\end{equation}
\end{defi}
An immediate consequence of this definition is the following
\begin{prop}\label{PAugmentedRackBialgebraIsRackBialgebra}
 Let $(B,\Phi,H,\ell)$ be an augmented rack bialgebra. Then the 
 $C^3$-coalgebra $(B,\epsilon,\mathbf{1})$ will become a 
 rack bialgebra
by means of the multiplication
\begin{equation}\label{EqDefRackBialgMultiplicationFromAugmented}
    a\triangleright b := \Phi(a).b
\end{equation}
for all $a,b\in B$. In particular, each Hopf algebra $H$ becomes
an augmented rack bialgebra via $(H,\mathrm{id}_H,H,\mathrm{ad})$. In general,
for each augmented rack bialgebra the map
$\Phi:B\to H$ is a morphism of rack bialgebras.
\end{prop}

\begin{prf} 
 We check first that $\triangleright$ is a morphism of $C^3$-coalgebras
 $B\otimes B\to B$: Let $a,b\in B$, then --thanks to the fact that
 the action $\ell$ and the maps $\Phi$ are coalgebra morphisms--
 \begin{eqnarray*}
  \Delta\big(\mu(a\otimes b)\big)
   & = & \Delta(a\triangleright b)
    = \Delta\big(\Phi(a).b\big)
    = \sum_{(\Phi(a)),(b)} 
      \Big(\big(\Phi(a)^{(1)}\big).b^{(1)}\Big)\otimes
      \Big(\big(\Phi(a)^{(2)}\big).b^{(2)}\Big) \\
      & = &\sum_{(a),(b)} 
      \Big(\big(\Phi(a^{(1)})\big).b^{(1)}\Big)\otimes
      \Big(\big(\Phi(a^{(2)})\big).b^{(2)}\Big) \\
      & = & \sum_{(a),(b)} \big(a^{(1)}\triangleright b^{(1)}\big)\otimes
      \big(a^{(2)}\triangleright b^{(2)}\big)
 \end{eqnarray*}
 whence $\mu$ is a morphism of coalgebras. Clearly
 \[
      \epsilon(a\triangleright b)
          = \epsilon\big(\Phi(a).b\big)=\epsilon_H\big(\Phi(a)\big)\epsilon(b)
          = \epsilon(a)\epsilon(b)
 \]
 whence $\mu$ preserves counits.\\
 We shall next compute both sides of the self-distributivity identity
 (\ref{EqDefAlgSelfDistributivity}) to get an idea: For all $a,b,c\in B$
 \[
    a\triangleright (b\triangleright c)
    =\Phi(a).\big(\Phi(b).c\big)
    =\big(\Phi(a)\Phi(b)\big).c,
 \]
 and
 \begin{eqnarray*}
   \sum_{(a)} (a^{(1)}\triangleright b)\triangleright
       (a^{(2)} \triangleright c) & = &
      \sum_{(a)}\big(\Phi(a^{(1)}).b\big) \triangleright
         \big(\Phi(a^{(2)}).c\big) \\
     & = &  \sum_{(a)}\Big(\Phi\big(\Phi(a^{(1)}).b\big)\Big).
         \big(\Phi(a^{(2)}).c\big) \\
      & = &  
      \sum_{(a)}
      \Big(\Phi\big(\Phi(a^{(1)}).b\big)\Phi(a^{(2)})\Big).c,      
 \end{eqnarray*}
and we compute, using the fact that $\Phi$ is a morphism of 
$C^3$-coalgebras,
\begin{eqnarray*}
  \sum_{(a)}
      \Phi\big(\Phi(a^{(1)}).b\big)\Phi(a^{(2)})
      & = &
      \sum_{(a)} 
      \Phi\big((\Phi(a)^{(1)}).b\big)(\Phi(a)^{(2)}) \\
      & 
      \stackrel{(\ref{EqDefAugBialgPhiIntertwinesActionWithAdjoint})}{=} 
          & 
       \sum_{(a)} 
       \Big(\mathrm{ad}_{\Phi(a)^{(1)}}\big(\Phi(b)\big)\Big)
       (\Phi(a)^{(2)}) \\
       & = &
       \sum_{(a)} 
       \big(\Phi(a)^{(1)}\big)\Phi(b)\Big(S\big(\Phi(a)^{(2)}\big)\Big)
          (\Phi(a)^{(3)}) \\
         & = &
       \sum_{(a)} 
       \big(\Phi(a)^{(1)}\big)\Phi(b)\mathbf{1}_H
           \epsilon_H\big(\Phi(a)^{(2)}\big)
        \\  
        & =& \Phi(a)\Phi(b),
\end{eqnarray*}
which proves the self-distributivity identity. Moreover we have
\[
      \mathbf{1}_B\triangleright a  = 
        \Phi(\mathbf{1}).a = \mathbf{1}_H.a = a,
\]
and
\[
      a \triangleright \mathbf{1}
      = \Phi(a).\mathbf{1} 
      \stackrel{(\ref{EqDefAugBialgHPreservesE})}{=}
      \epsilon_H\big(\Phi(a)\big)\mathbf{1}
      =\epsilon_B(a)\mathbf{1},
\]
whence the $C^3$-coalgebra becomes a rack bialgebra. 
\end{prf}

\begin{exem}  \label{augmented_rack_algebra}
Exactly in the same way as a pointed rack gives rise to a rack bialgebra $K[X]$,
an augmented pointed rack $p:X\to G$ gives rise to an augmented rack bialgebra
$p:K[X]\to K[G]$. 
\end{exem} 

\begin{rem}
Motivated by the fact that the augmented racks $p:X\to G$ are exactly the Yetter-Drinfeld modules
over the (set-theoretical) Hopf algebra $G$, we may ask whether augmented rack bialgebras are 
Yetter-Drinfeld modules.
 
In fact, any cocommutative augmented rack bialgebra $(B,\Phi,H,\ell)$ gives rise 
to a Yetter-Drinfeld module over the Hopf algebra $H$. Indeed, $B$ is a left $H$-module
via $\ell$, and becomes a left $H$-comodule via
$$\rho:B\stackrel{\triangle_B}{\to}B\otimes B\stackrel{\Phi\otimes\id_B}{\to}H\otimes B.$$
Now, in Sweedler notation, the coaction is denoted for all $b\in B$ by
$$\rho(b)\,=\,\sum_{(b)}b_{(-1)}\otimes b_{(0)}\,\in\,H\otimes B.$$
Then the {\it Yetter-Drinfeld compatibility relation} reads
$$\sum_{(h.b)}(h.b)_{(-1)}\otimes(h.b)_{(0)}\,=\,\sum_{(b),(h)}h^{(1)}b_{(-1)}S(h^{(3)})\otimes h^{(2)}.b_{(0)}.$$
This relation is true in our case, because $\ell$ is a morphism of coalgebras and 
is sent to the adjoint action via $\Phi$. 

Conversely, given a Yetter-Drinfeld module $C$ over a Hopf algebra $H$,
together with a linear form $\epsilon_C:C\to K$
satisfying $\epsilon_C(h.c) =\epsilon_H(h)\epsilon_C(c)$, then define a map 
$\Phi:C\to H$ by
$$\Phi\,:=\,(\id_H\otimes\epsilon_C)\circ\rho.$$ 
The map $\Phi$ intertwines the left action on $C$ and the adjoint action on $H$ 
thanks to the Yetter-Drinfeld condition. 

Now define a rack product for all $x,y\in C$ by 
$$x\rhd y\,=\,\Phi(x).y,$$ 
then we obtain a Yetter-Drinfeld version of self-distributivity
$$x\rhd(y\rhd z)\,=\,\sum_{(x)}(x_{(-1)}.y)\rhd(x_{(0)}\rhd z),$$ 
as there is no comultiplication on $C$. 

The fact that $\Phi$ is a morphism of coalgebras is then replaced by the identity
$$(\id_H\otimes \Phi)\circ\rho\,=\,\Delta_H\circ\Phi,$$
which one needs to demand.  

Finally, one needs a unit $1_C\in C$ such that for all $h\in H$,
$h.1_C=\epsilon_H(h)1_C$, $\epsilon_C(1_C)=1_K$, 
$\rho(1_C)=1_H\otimes 1_C$, and $\Phi(1_C)= 1_H$.  
This is somehow the closest one can get to a rack bialgebra without
having a compatible $C^3$ coalgebra structure on $C$.  
\end{rem} 
 
\subsection{(Augmented) rack bialgebras for any Leibniz algebra}

In this subsection, we will suppose $\Q\subset K$. 
Let $(\mathfrak{h}, [~,~])$ be a \emph{Leibniz algebra over $K$}, i.e.
$\mathfrak{h}$ is a $K$-module equipped with a $K$-linear map 
$[~,~]:\mathfrak{h}\otimes \mathfrak{h}\to\mathfrak{h}$ satisfying
the (left) Leibniz identity (\ref{LeibnizIdentity}).

Recall first that each Lie algebra over $K$ is a Leibniz algebra
giving rise to a functor from the category of all Lie algebras to the 
category of all Leibniz algebras.

Furthermore, recall that each Leibniz algebra has two canonical 
$K$-submod\-ules
\begin{eqnarray}
   Q(\mathfrak{h})&:=&
   \big\{x\in\mathfrak{h}~|~\exists~N\in\mathbb{N}\setminus\{0\},~
    \exists~\lambda_1,\ldots,\lambda_N\in K,~\exists~x_1,\ldots,x_N~
    \nonumber\\
    & &
   ~~~~~~~~~~~~~~~~~~ 
       \mathrm{such~that~}x=\sum_{r=1}^N\lambda_r[x_r,x_r]\big\}, \\
    \mathfrak{z}(\mathfrak{h})
    & := & \big\{x\in\mathfrak{h}~|~\forall~y\in\mathfrak{h}
                     :~[x,y]=0\big\}.
\end{eqnarray}
It is well-known and not hard to deduce from the Leibniz identity 
that both
$Q(\mathfrak{h})$ and $\mathfrak{z}(\mathfrak{h})$ are two-sided
abelian ideals of $(\mathfrak{h}, [~,~])$, that
$Q(\mathfrak{h})\subset \mathfrak{z}(\mathfrak{h})$, and that
the quotient Leibniz algebras
\begin{equation}
    \overline{\mathfrak{h}}:=\mathfrak{h}/Q(\mathfrak{h})
    ~~~\mathrm{and}~~~
       \mathfrak{g}(\mathfrak{h}):=
       \mathfrak{h}/\mathfrak{z}(\mathfrak{h})
\end{equation}
are Lie algebras. Since the ideal $Q(\mathfrak{h})$ is clearly
mapped
into the ideal $Q(\mathfrak{h}')$ by any morphism of Leibniz algebras
$\mathfrak{h}\to\mathfrak{h}'$ (which is a priori not the case
for $\mathfrak{z}(\mathfrak{h})$ !), there is an obvious functor
$\mathfrak{h}\to \overline{\mathfrak{h}}$ from the category of
all Leibniz algebras to the category of all Lie algebras.

In order to perform the following constructions of rack bialgebras for
any given Leibniz
algebra $(\mathfrak{h},[~,~])$, choose first a two-sided ideal
$\mathfrak{z}\subset \mathfrak{h}$ such that
\begin{equation}
    Q(\mathfrak{h})\subset \mathfrak{z} \subset 
    \mathfrak{z}(\mathfrak{h}),
\end{equation}
let $\mathfrak{g}$ denote the quotient Lie algebra 
$\mathfrak{h}/\mathfrak{z}$, and let $p:\mathfrak{h}\to \mathfrak{g}$
be the natural projection. The data of $\mathfrak{z}\subset\mathfrak{h}$,
i.e. of a Leibniz algebra $\mathfrak{h}$ together with an ideal $\mathfrak{z}$
such that $Q(\mathfrak{h})\subset\mathfrak{z}\subset\mathfrak{z}(\mathfrak{h})$,
could be called an {\it augmented Leibniz algebra}. Thus we are actually associating 
an augmented rack bialgebra to every augmented Leibniz algebra. 
In fact, we will see that this augmented rack bialgebra does not depend on the choice of the ideal 
$\mathfrak{z}$ and therefore refrain from introducing augmented Leibniz algebras in a more 
formal way.

The Lie algebra $\mathfrak{g}$ naturally
acts as derivations on $\mathfrak{h}$ by means of 
(for all $x,y\in\mathfrak{h}$)
\begin{equation} \label{EqDefLeibnizHIsAGModule}
    p(x).y := [x,y]=:\mathrm{ad}_x(y)
\end{equation}
because $\mathfrak{z} \subset \mathfrak{z}(\mathfrak{h})$. Note that
\begin{equation}\label{EqCompHModLeftCenterCongAd}
 \mathfrak{h}/\mathfrak{z}(\mathfrak{h})\cong 
   \big\{\mathrm{ad}_x\in\mathrm{Hom}_K(\mathfrak{h},\mathfrak{h})~|~
         x\in\mathfrak{h}\big\}.
\end{equation}
as Lie algebras.

Consider 
now the $C^5$-coalgebra 
$(B=\mathsf{S}(\mathfrak{h}),\Delta,\epsilon,\mathbf{1})$.
Here $\mathsf{S}(\mathfrak{h})$ is the symmetric algebra and coalgebra on 
the vector space $\mathfrak{h}$, and $\mathsf{S}(\mathfrak{h})$ is
actually a commutative cocommutative Hopf algebra over $K$ with respect
to the symmetric multiplication $\bullet$. The linear map 
$p:\mathfrak{h}\to \mathfrak{g}$ induces a unique morphism of
Hopf algebras
\begin{equation}
   \tilde{\Phi}=\mathsf{S}(p):
     \mathsf{S}(\mathfrak{h})\to \mathsf{S}(\mathfrak{g})
\end{equation}
satisfying
\begin{equation}\label{EqDefMapTildePhi}
    \tilde{\Phi}(x_1\bullet\cdots\bullet x_k)=p(x_1)\bullet\cdots\bullet p(x_k)
\end{equation}
for any nonnegative integer $k$ and $x_1,\ldots,x_k\in \mathfrak{h}$.
In other words, the association $\mathsf{S}:V\to \mathsf{S}(V)$
is a functor from the category of all $K$-modules to the category
of all commutative unital $C^5$-coalgebras.
Consider now the universal enveloping algebra $\mathsf{U}(\mathfrak{g})$
of the Lie algebra $\mathfrak{g}$. Since $\mathbb{Q}\subset K$ by 
assumption, the Poincar\'{e}-Birkhoff-Witt Theorem (in short: PBW) holds
(see e.g. \cite[Appendix]{Qui69}). More precisely, the symmetrisation
map $\omega:\mathsf{S}(\mathfrak{g})\to \mathsf{U}(\mathfrak{g})$,
defined by
\begin{equation}\label{EqDefSymmetrisationMap}
   \omega(\mathbf{1}_{\mathsf{S}(\mathfrak{g})})
      = \mathbf{1}_{\mathsf{U}(\mathfrak{g})},~~~\mathrm{and}~~~
    \omega(\xi_1\bullet\cdots\bullet \xi_k)
    = \frac{1}{k!}\sum_{\sigma\in S_k}
        \xi_{\sigma(1)}\cdots\xi_{\sigma(k)},
\end{equation}
see e.g. \cite[p.80, eqn (3)]{Dix74}, is an isomorphism of $C^5$-coalgebras 
(in general not of associative algebras). We now need an action of the 
Hopf algebra $H=\mathsf{U}(\mathfrak{g})$ on $B$, and an intertwining
map $\Phi:B\to\mathsf{U}(\mathfrak{g})$. In order to get this, we first look
at $\mathfrak{g}$-modules: The $K$-module $\mathfrak{h}$ is a 
$\mathfrak{g}$-module by means of eqn (\ref{EqDefLeibnizHIsAGModule}), the 
Lie algebra $\mathfrak{g}$ is a $\mathfrak{g}$-module via its adjoint
representation, and the linear map $p:\mathfrak{h}\to\mathfrak{g}$ is a 
morphism of $\mathfrak{g}$-modules since $p$ is a morphism of Leibniz algebras.
Now $\mathsf{S}(\mathfrak{h})$ and $\mathsf{S}(\mathfrak{g})$ are
$\mathfrak{g}$-modules in the usual way, i.e. for all 
$k\in\mathbb{N}\setminus \{0\}$, $\xi,\xi_1,\ldots,\xi_k\in\mathfrak{g}$, and
$x_1\ldots,x_k\in\mathfrak{h}$
\begin{eqnarray}
 \xi.(x_1\bullet\cdots\bullet x_k)
 & := & \sum_{r=1}^k
  x_1\bullet\cdots\bullet (\xi.x_r)\bullet \cdots \bullet x_k, 
     \label{EqDefGActionOnSH}\\
  \xi.(\xi_1\bullet\cdots\bullet \xi_k)
 & := & \sum_{r=1}^k
  \xi_1\bullet\cdots\bullet [\xi.\xi_r]\bullet \cdots \bullet \xi_k,
  \label{EqDefGActionOnSG}
\end{eqnarray}
and of course $\xi.\mathbf{1}_{\mathsf{S}(\mathfrak{h})}=0$ and
$\xi.\mathbf{1}_{\mathsf{S}(\mathfrak{g})}=0$. Recall that 
$\mathsf{U}(\mathfrak{g})$ is a $\mathfrak{g}$-module via the adjoint 
representation $\mathrm{ad}_\xi(u)=\xi.u= \xi u - u\xi$ 
(for all $\xi\in\mathfrak{g}$ and
all $u\in \mathsf{U}(\mathfrak{g})$).\\
It is easy to see that the map $\tilde{\Phi}$ (\ref{EqDefMapTildePhi}) is a 
morphism
of $\mathfrak{g}$-modules, and it is well-known that the symmetrization map 
$\omega$ (\ref{EqDefSymmetrisationMap}) is also a morphism of
$\mathfrak{g}$-modules, see e.g. \cite[p.82, Prop. 2.4.10]{Dix74}.
Define the $K$-linear map $\Phi:\mathsf{S}(\mathfrak{h})\to 
\mathsf{U}(\mathfrak{g})$ by the composition
\begin{equation}\label{EqDefPBWCircSp}
    \Phi:= \omega\circ \tilde{\Phi}.
\end{equation}
Then $\Phi$ is a map of $C^5$-coalgebras and a map of $\mathfrak{g}$-modules.
Thanks to the universal property of the universal enveloping algebra, it 
follows
that $\mathsf{S}(\mathfrak{h})$ and $\mathsf{U}(\mathfrak{g})$ are 
left $\mathsf{U}(\mathfrak{g})$-modules, via
(for all $\xi_1,\ldots,\xi_k\in\mathfrak{g}$, and for all $a\in 
\mathsf{S}(\mathfrak{h})$) 
\begin{equation}\label{EqDefUGActionOnSH}
  (\xi_1\cdots\xi_k).a = \xi_1.(\xi_2.(\cdots \xi_k.a)\cdots)
\end{equation}
and the usual adjoint representation (\ref{EqSecondDefAdjointRep})
(for all $u\in \mathsf{U}(\mathfrak{g})$) 
\begin{equation}
 \mathrm{ad}_{\xi_1\cdots\xi_k}(u)
 =\big(\mathrm{ad}_{\xi_1}\circ\cdots\circ\mathrm{ad}_{\xi_k}\big)(u),
\end{equation}
and that $\Phi$ intertwines the $\mathsf{U}(\mathfrak{g})$-action on
$C=\mathsf{S}(\mathfrak{h})$ with the adjoint action of 
$\mathsf{U}(\mathfrak{g})$ on itself.\\
Finally it is a routine check using the above identities 
(\ref{EqDefGActionOnSH}) and (\ref{EqSecondDefAdjointRep}) that 
$\mathsf{S}(\mathfrak{h})$ becomes a module coalgebra.\\
We can resume the preceding considerations in the following
\begin{theo}\label{TLeibnizToUARInfinityOfLeibniz}
  Let $(\mathfrak{h},[~,~])$ be a Leibniz algebra over $K$, let
  $\mathfrak{z}$ be a two-sided ideal of $\mathfrak{h}$ such that
  $Q(\mathfrak{h})\subset \mathfrak{z} \subset \mathfrak{z}(\mathfrak{h})$,
  let $\mathfrak{g}$ denote the quotient Lie algebra 
  $\mathfrak{h}/\mathfrak{z}$ by $\mathfrak{g}$, and let 
  $p:\mathfrak{h}\to \mathfrak{g}$ be the canonical projection.
  \begin{enumerate}
  \item Then there is a canonical $\mathsf{U}(\mathfrak{g})$-action $\ell$ on
  the $C^5$-coalgebra $B:=\mathsf{S}(\mathfrak{h})$ (making it into a module
  coalgebra leaving invariant $\mathbf{1}$) and a canonical lift of
  $p$ to a map of $C^5$-coalgebras, $\Phi:\mathsf{S}(\mathfrak{h})\to
  \mathsf{U}(\mathfrak{g})$ such that eqn 
  (\ref{EqDefAugBialgPhiIntertwinesActionWithAdjoint}) holds.\\
  Hence the quadruple $(\mathsf{S}(\mathfrak{h}),\Phi,
  \mathsf{U}(\mathfrak{g}), \ell)$ is an augmented rack bialgebra 
  whose associated Leibniz algebra is equal to $(\mathfrak{h},[~,~])$
  (independently of the choice of $\mathfrak{z}$).\\
  The resulting rack multiplication $\mu$ of $\mathsf{S}(\mathfrak{h})$ 
  (written
  $\mu(a\otimes b)=a\triangleright b$) is also independent on the choice of
  $\mathfrak{z}$ and is explicitly given
  as follows for all positive integers $k,l$ and
  $x_1,\ldots,x_k$, $y_1,\ldots,y_l\in\mathfrak{h}$:
  \begin{equation}   \label{adjoint_action_s}
      \big(x_1\bullet\cdots\bullet x_k)\triangleright
      \big(y_1\bullet\cdots\bullet y_l)
      =\frac{1}{k!}\sum_{\sigma\in S_k}
        \big(\mathrm{ad}^s_{x_{\sigma(1)}}\circ\cdots\circ 
          \mathrm{ad}^s_{x_{\sigma(k)}}\big)
          \big(y_1\bullet\cdots\bullet y_l)
  \end{equation}
  where $\mathrm{ad}^s_x$ denotes the action of the Lie algebra
  $\mathfrak{h}/\mathfrak{z}(\mathfrak{h})$ (see eqn 
  (\ref{EqCompHModLeftCenterCongAd}))
  on $\mathsf{S}(\mathfrak{h})$ according to eqn (\ref{EqDefGActionOnSH}).

  \item In case $\mathfrak{z}=Q(\mathfrak{h})$, the construction mentioned
  in {\it 1.} is
  a functor $\mathfrak{h}\to\mathsf{UAR}^\infty(\mathfrak{h})$ from the 
  category of all Leibniz algebras to the 
  category of all augmented rack bialgebras associating to ${\mathfrak h}$
  the rack bialgebra 
$$\mathsf{UAR}^\infty(\mathfrak{h}):=(\mathsf{S}(\mathfrak{h}),\Phi,
  \mathsf{U}(\mathfrak{g}),\ell)$$ 
  and to each morphism
  $f$ of Leibniz algebras the pair $\big(\mathsf{S}(f),
  \mathsf{U}(\overline{f})\big)$ where $\overline{f}$ is the induced Lie 
  algebra morphism.
  \item For each nonnegative integer $k$, the above construction restricts
  to each subcoalgebra of order $k$, $\mathsf{S}(\mathfrak{h})_{(k)}
   =\oplus_{r=0}^k \mathsf{S}^r(\mathfrak{h})$, to define an augmented
   rack bialgebra  $(\mathsf{S}(\mathfrak{h})_{(k)},\Phi_{(k)},
  \mathsf{U}(\mathfrak{g}), 
  \ell|_{\mathsf{U}(\mathfrak{g})\otimes \mathsf{S}(\mathfrak{h})_{(k)}})$
  which in case  $\mathfrak{z}=Q(\mathfrak{h})$ defines a functor
  $\mathfrak{h}\to
  \mathsf{UAR}_{(k)}(\mathfrak{h}):=
  \big(\mathsf{UAR}^\infty(\mathfrak{h})\big)_{(k)}$ from the 
  category of all Leibniz algebras to the 
  category of all augmented rack bialgebras.
  \end{enumerate}
\end{theo}
\begin{prf}
  {\it 1.} All the statements except the last two ones
  have already been proven. Note that
    for all $x,y\in\mathfrak{h}$ we have by definition
     \[
         [x,y]= p(x).y =x\triangleright y,
     \]
     independently of the chosen ideal $\mathfrak{z}$. Moreover we compute
     \begin{eqnarray*}
       \lefteqn{\big(x_1\bullet\cdots\bullet x_k)\triangleright
      \big(y_1\bullet\cdots\bullet y_l)}\\
      & = & \big((\omega\circ \tilde{\Phi})(x_1\bullet\cdots\bullet x_k)\big).
          \big(y_1\bullet\cdots\bullet y_l)\\
      & = & \frac{1}{k!}\sum_{\sigma\in S_k}
        \big(p(x_{\sigma(1)})\cdots
          p(x_{\sigma(k)})\big)
          .\big(y_1\bullet\cdots\bullet y_l),    
     \end{eqnarray*}
    which gives the desired formula since for all $x\in \mathfrak{h}$ and
    $a\in\mathsf{S}(\mathfrak{h})$, we have
    \[
        p(x).a=\mathrm{ad}^s_x(a).
    \]
  {\it 2.} Let $f:\mathfrak{h}\to \mathfrak{h}'$ be a morphism
  of Leibniz algebras, and let $\overline{f}:
  \overline{\mathfrak{h}}\to \overline{\mathfrak{h}'}$ be the induced 
  morphism of Lie algebras. Hence we get
  \begin{equation}\label{EqCompPPrimeCirclePsiEqualsOverlinePsiCircleP}
       p'\circ f = \overline{f}\circ p
  \end{equation}
  where $p':\mathfrak{h}'\to \overline{\mathfrak{h}'}$ denotes the 
  corresponding projection modulo $Q(\mathfrak{h}')$.
  Let $\mathsf{S}(f):\mathsf{S}(\mathfrak{h})\to 
  \mathsf{S}(\mathfrak{h}')$,
  $\mathsf{S}(\overline{f}):
  \mathsf{S}(\overline{\mathfrak{h}})\to 
  \mathsf{S}(\overline{\mathfrak{h}'})$, 
  and
  $\mathsf{U}(\overline{f}):
  \mathsf{U}(\overline{\mathfrak{h}})\to 
  \mathsf{U}(\overline{\mathfrak{h}'})$ 
  the induced maps of Hopf algebras, i.e. $\mathsf{S}(f)$ (resp. 
  $\mathsf{S}(\overline{f})$) satisfies eqn 
  (\ref{EqDefMapTildePhi}) (with $p$ replaced by $f$ (resp. 
  by $\overline{f}$)), and $\mathsf{U}(\overline{f})$ satisfies
  \[
    \mathsf{U}(\overline{f}) \big(\xi_1\cdots \xi_k\big)=
        \overline{f}(\xi_1)\cdots \overline{f}(\xi_k)
  \]
  for all positive integers $k$ and $\xi_1,\ldots,\xi_k\in
  \overline{\mathfrak{h}}$. If $\omega:
  \mathsf{S}(\overline{\mathfrak{h}})\to
  \mathsf{U}(\overline{\mathfrak{h}})$
  and $\omega':
  \mathsf{S}(\overline{\mathfrak{h}'})\to
  \mathsf{U}(\overline{\mathfrak{h}'})$ denote the corresponding 
  symmetrisation 
  maps (\ref{EqDefSymmetrisationMap}) then it is easy to see from the 
  definitions that
  \[
      \omega'\circ \mathsf{S}(\overline{f})
        = \mathsf{U}(\overline{f})\circ \omega.
  \]
  Equation (\ref{EqCompPPrimeCirclePsiEqualsOverlinePsiCircleP}) implies
  \[
     \tilde{\Phi}'\circ \mathsf{S}(f) = \mathsf{S}(p')\circ \mathsf{S}(f)
       = \mathsf{S}(\overline{f})\circ  \mathsf{S}(p)
            =\mathsf{S}(\overline{f})\circ\tilde{\Phi},
  \]
  and composing from the left with $\omega'$ yields the equation
  \begin{equation}\label{EqCompPhiPrimePsiEqualsOverlinePsiPhi}
   \Phi'\circ \mathsf{S}(f) = \mathsf{U}(\overline{f})\circ  \Phi.
  \end{equation}
  Moreover for all $x,y\in\mathfrak{h}$ we have, since
  $f$ is a morphism of Leibniz algebras,
  \[
     f\big(p(x).y\big) =
       f\big([x,y]\big) = \big[f(x),f(y)\big]'
        = p'\big(f(x)\big).f(y)
        = \overline{f}\big(p(x)\big).f(y),
  \]
  hence for all $\xi\in\overline{h}$
  \[
     f(\xi.y)=\big(\overline{f}(\xi)\big).\big(f(y)\big),
  \]
  and upon using eqn (\ref{EqDefGActionOnSH}) we get for all $a\in
  \mathsf{S}(\mathfrak{h})$
  \[
      \mathsf{S}(f)\big(\xi.a\big)
      =\big(\overline{f}(\xi)\big).\Big(\mathsf{S}(f)\big(a\big)\Big), 
  \]
   showing finally
  for all 
  $u\in\mathsf{U}(\mathfrak{h})$ and all $a\in\mathsf{S}(\mathfrak{h})$
  \begin{equation}\label{EqCompIntertwiningOfModuleActions}
   \mathsf{S}(f)\big(u.a\big) 
   = \big(\mathsf{U}(\overline{f})(u)\big).\Big(\mathsf{S}(f)\big(a\big)\Big).
  \end{equation}
  Associating to every Leibniz algebra
  $(\mathfrak{h}, [~,~])$ the above defined augmented rack bialgebra
  $(\mathsf{S}(\mathfrak{h}),\Phi,\mathsf{U}(\overline{\mathfrak{h}}),
   \ell)$, and to every morphism $\psi:\mathfrak{h}\to \mathfrak{h}'$
   of Leibniz algebras the pair of $K$-linear maps 
   $\big(\Psi=\mathsf{S}(\psi),\overline{\Psi}=
   \mathsf{U}(\overline{\psi})\big)$, we can easily check that 
   $\Psi$ is a morphism
   of $C^5$-coalgebras, $\overline{\Psi}$ is a morphism of Hopf algebras,
   such that the two relevant diagrams 
  (\ref{EqDefMorphAugmentedRackBialg}) commute which easily follows from
   (\ref{EqCompPhiPrimePsiEqualsOverlinePsiPhi}) and
   (\ref{EqCompIntertwiningOfModuleActions}). The rest of the functorial 
   properties is a routine check.\\  
   {\it 3.} By definition, the $\mathsf{U}(\mathfrak{g})$-action on
   $\mathsf{S}(\mathfrak{h})$ (cf. eqs (\ref{EqDefGActionOnSH})
   and (\ref{EqDefUGActionOnSH})) leaves invariant each $K$-submodule
   $\mathsf{S}^r(\mathfrak{h})$ for each nonnegative integer $r$ whence it leaves
   invariant each subcoalgebra of order $k$, $\mathsf{S}(\mathfrak{h})_{(k)}$.
   It follows that the construction restricts well.
\end{prf}
 
\begin{defi} \label{definition_UARinfini}
The rack bialgebra $\mathsf{UAR}^\infty(\mathfrak{h})$ which is by the above 
theorem canonically associated to each Leibniz algebra $\mathfrak{h}$
is called the universal augmented rack bialgebra. 
\end{defi}

\begin{rem}
This theorem should be compared to Proposition 3.5 in \cite{CCES08}. 
In \cite{CCES08},
the authors work with the vector space $N:=K\oplus{\mathfrak h}$, while we 
work with the whole symmetric algebra on the Leibniz algebra. In some sense, 
we extend their Proposition 3.5 ``to all orders'', hence the name 
$\mathsf{UAR}^\infty(\mathfrak{h})$. The subset $N=K\oplus{\mathfrak h}$ becomes
a sub rack bialgebra denoted by $\mathsf{UAR}(\mathfrak{h})$. It turns out that 
$\mathsf{UAR}(\mathfrak{h})$ is already enough to obtain a 
left-adjoint to the functor of primitives and hence universality, 
see \cite{ABRW}. 

The above rack bialgebra $\mathsf{UAR}^\infty(\mathfrak{h})$ associated to a 
Leibniz algebra ${\mathfrak h}$ can be seen
as one version of an {\it enveloping algebra} of 
the Leibniz algebra ${\mathfrak h}$.
The link to the \emph{universal enveloping algebra of
$\mathfrak{h}$}, $\mathfrak{h}\otimes \mathsf{U}(\mathfrak{\overline{h}})$, as a dialgebra
(in the sense of 
of Loday-Pirashvili) has been elucidated in \cite{ABRW}. 
\end{rem}

We shall close the subsection with a geometric explanation of some of the 
structures appearing here: Let
$\big(\mathfrak{h},[~,~]\big)$ be a real finite-dimensional Leibniz algebra.
Then for any real number $\hbar$, there is the following Lie rack structure on 
the manifold $\mathfrak{h}$
defined by
\begin{equation}\label{EqDefStandardLieRackStructureOnLeibniz}
    x\blacktriangleright_\hbar y := e^{\hbar\mathrm{ad}_x}(y)
\end{equation}
For later use we note that on the space $\mathfrak{h}[[\hbar]]$ of all formal
power series the above formula makes sense if $x,y$ are also formal
power series.\\
Moreover, pick a two-sided ideal $\mathfrak{z}\subset \mathfrak{h}$ with
$Q(\mathfrak{h})\subset \mathfrak{z}\subset \mathfrak{z}(\mathfrak{h})$ so
that the quotient algebra $\mathfrak{g}:=\mathfrak{h}/\mathfrak{z}$ is a Lie 
algebra. Let $p:\mathfrak{h}\to \mathfrak{g}$ be the canonical projection.
Let $G$ be the connected simply connnected Lie group having Lie algebra
$\mathfrak{g}$. Since $\mathfrak{g}$ acts on $\mathfrak{h}$ as derivations,
there is a unique Lie group action $\ell$ of $G$ on $\mathfrak{h}$ by automorphisms
of Leibniz algebras. Consider the smooth map 
\begin{equation}\label{EqDefPhiForStandardAugmentedLieRack}
    \phi: \mathfrak{h}\to G: x\mapsto \exp\big(p(x)\big).
\end{equation}
Clearly $\phi(g.x)=g\phi(x)g^{-1}$ for all $x\in\mathfrak{h}$ and $g\in G$
whence $(\mathfrak{h},\phi,G,\ell)$ is an augmented Lie rack, and it is not hard
to see that the Lie rack structure coincides with 
(\ref{EqDefStandardLieRackStructureOnLeibniz}) for $\hbar=1$.
The following theorem is shown in \cite{ABRW}:
\begin{theo}
 The $C^5$-rack bialgebra associated to the augmented Lie rack
 $(\mathfrak{h},\phi,G,\ell)$ by means of the functor
 $F$, described in Example 1.4,
 is isomorphic to the universal envelopping algebra of infinite order,
 $\mathsf{UAR}^\infty(\mathfrak{h})$ 
 (Definition \ref{definition_UARinfini}).
\end{theo}

\subsection{Quantum racks}

In the article \cite{DW2013}, the authors are interested in the generalized 
Poisson 
manifold given by the linear dual ${\mathfrak h}^*$ of a finite-dimensional 
Leibniz algebra ${\mathfrak h}$. Recall that the dual ${\mathfrak g}^*$
of a Lie algebra ${\mathfrak g}$ is a Poisson manifold with the 
Kostant-Kirillov-Souriau bracket, given for 
$f,g\in{\mathcal C}^{\infty}({\mathfrak g}^*)$ by
$$\{f,g\}(\xi):=\langle \xi,[df,dg]\rangle\,=\,\sum_{i,j,k}c_{ij}^k
\frac{\partial f}{\partial x_i}(\xi)\frac{\partial g}{\partial x_j}(\xi)X_k,$$
where $df$ and $dg$ are seen as linear functions on ${\mathfrak g}^*$, 
i.e. elements of ${\mathfrak g}$. The $c_{ij}^k$ are the structure constants
of the Lie bracket with respect to a certain basis of $(X_k)$. 
In the same manner, the dual of a Leibniz algebra ${\mathfrak h}^*$ carries 
a bracket given for $f,g\in{\mathcal C}^{\infty}({\mathfrak h}^*)$ by
$$\{f,g\}(\xi):=-\langle \xi,[df(0),dg]\rangle\,=\,-\sum_{i,j,k}
c_{ij}^k\frac{\partial f}{\partial x_i}(0)\frac{\partial g}{\partial x_j}
(\xi)X_k,$$
i.e. with respect to the above formula, the bracket is here evaluated in 
$0\in{\mathfrak h}^*$. The bracket is neither antisymmetric, nor a 
biderivation, nor does it satisfy some Jacobi/Leibniz identity in general. 
Nevertheless, \cite{DW2013} shows that this is the natural bracket
one encounters when following the standard deformation quantization procedure
for the dual of Lie algebras.   

Going into details, a finite dimensional Lie algebra ${\mathfrak g}$ admits a 
local Lie group $G$ giving rise to a (local) multiplication
$\mu:G\times G\to G$. The cotangent lift $T^*\mu$ can be interpreted as a 
Lagrangian submanifold of the triple product
$\overline{T^*G}\times\overline{T^*G}\times T^*G$. It constitutes thus 
a {\it symplectic micromorphism} between germs of symplectic manifolds, see
\cite{CDW}. When such a symplectic micromorphism is given by a 
{\it generating function} $S$, the function $S$ can be taken as an ``action'' 
in some oscillatory integral giving the deformation quantization of the 
corresponding Poisson manifold by Fourier Integral Operators.
In \cite{DW2013}, this quantization scheme is adapted to 
the dual of Leibniz algebras. The key ingredient is the integration
of ${\mathfrak h}$ into a Lie rack where for $x,y\in{\mathfrak h}$, the 
rack product is given by 
$$X\blacktriangleright Y\,=\,e^{{\rm ad}_X}(Y).$$
We take then $S_{\rhd}(X,Y,\xi)=\langle \xi,
e^{{\rm ad}_X}(Y)\rangle$ as a generating function, and the oscillatory 
integral is written 
$$Q_{\rhd}(f\otimes g)(\xi)\,=\,\int_{{\mathfrak{g}}\times{\mathfrak{g}}}
\widehat{f}
(X)\widehat{g}(Y)e^{\frac{i}{\hbar}S_{\rhd}(X,Y,\xi)}
\frac{dXdY}{(2\pi\hbar)^{n}}.$$
Here $n$ is the dimension of ${\mathfrak h}$ and $\widehat{f}$, 
$\widehat{g}$ are the asymptotic Fourier transforms of $f$ and $g$. 
The stationary phase series expansion of this integral gives then the 
corresponding star product. Its first term is the above generalized Poisson 
bracket.

The main theorem of \cite{DW2013} reads:

\begin{theo}  \label{theorem_geometric_defoquant}
The operation
\[
\rhd_{\hbar}:{\mathcal{C}}^{\infty}({\mathfrak{h}}^{*})[[\epsilon]]\otimes{\mathcal{C}}^{\infty}({\mathfrak{h}}^{*})[[\epsilon]]\to{\mathcal{C}}^{\infty}({\mathfrak{h}}^{*})[[\epsilon]]
\]
defined by
\[
f\rhd_{\hbar}g\,:=\, Q_{\rhd}(f\otimes g)
\]
is a quantum rack, i.e. 

(1) $\rhd_{\hbar}$ restricted to $U_{\mathfrak{h}}=\{E_{X}:=e^{\frac{i}{\hbar}X}\,|\, X\in{\mathfrak{h}}\}$
is a rack structure and
\[
e^{\frac{i}{\hbar}X}\rhd_{\hbar}e^{\frac{i}{\hbar}Y}\,=\, e^{\frac{i}{\hbar}{\rm conj}_{*}(X,Y)},
\]

(2) $\rhd_{\hbar}$ restricted to $\rhd_{\hbar}:U_{\mathfrak{h}}\times{\mathcal{C}}^{\infty}({\mathfrak{h}}^{*})\to{\mathcal{C}}^{\infty}({\mathfrak{h}}^{*})$
is a rack action and 
\[
(e^{\frac{i}{\hbar}X}\rhd_{\hbar}f)(\xi)\,=\,({\rm Ad}_{-X}^{*}f)(\xi).
\]

\end{theo}

\section{Deformation quantization of rack bialgebras}
 
\subsection{Algebraic deformation quantization of Leibniz algebras}

In this subsection, $K$ denotes the field of real numbers $\R$ or 
the field of complex numbers $\C$.

Let $({\mathfrak h},[\,,\,])$ be a finite dimensional Leibniz algebra of 
dimension $n$, and denote by ${\mathfrak h}^*$ its linear dual. 
In order to make computations more elementary we shall use a fixed basis
$e_1,\ldots,e_n$ of ${\mathfrak h}$, but it is a routine check that
all the relevant formulas are invariant under a change of basis.
Let $e^1,\ldots,e^n$ be 
the corresponding dual basis of ${\mathfrak h}^*$, i.e. by definition
$$
e^i(e_j)\,=\,\delta^i_{j},
$$
for all $i,j=1,\ldots,n$. Furthermore, let $c_{ij}^k$ for $i,j,k=1,\ldots,n$
be the structure constants of the Leibniz algebra ${\mathfrak h}$ with 
respect to the basis $e_1,\ldots,e_n$, i.e.
$$
c^i_{jk}\,=\,e^i([e_j,e_k])
$$
for all $i,j,k=1,\ldots,n$. We will denote by $x,y,z,\ldots$ elements of 
${\mathfrak h}$, while $\alpha,\beta,\gamma,\ldots$ will denote elements of 
${\mathfrak h}^*$. Denote by $\alpha_1,\ldots,\alpha_n$ the 
coordinates of $\alpha\in{\mathfrak h}^*$ with respect to the basis 
$e^1,\ldots,e^n$. 
For all $x\in{\mathfrak h}$, denote by 
$\hat{x}\in{\mathcal C}^{\infty}({\mathfrak h}^*,K)$ the 
linear function given by 
$$
\hat{x}(\alpha)\,:=\,\alpha(x),
$$
for all $\alpha\in{\mathfrak h}^*$. In the same vein, let $e^{\hat{x}}$ be the 
exponential function given by 
$$
e^{\hat{x}}(\alpha)\,:=\,e^{\alpha(x)}\,=\,e^{\hat{x}(\alpha)},
$$
for all $\alpha\in{\mathfrak h}^*$. For all integers $i=1,\ldots,n$, define
a first order differential operator $\widetilde{\rm ad}_i$ on smooth functions  
$f:{\mathfrak h}^*\to K$ by
$$
(\widetilde{\rm ad}_i(f))(\alpha)\,:=\,\sum_{j,k=1}^n\alpha_k\,c_{ij}^k\,
\frac{\partial f}{\partial\alpha_j}(\alpha).
$$ 
The following star-product formula, where $\hbar$ is 
a formal parameter (which may be replaced by a real number in situations where 
the formula is convergent), will render ${\mathfrak h}^*$ a {\it quantum rack}
in the sense of Theorem \ref{theorem_geometric_defoquant}, see also
\cite{DW2013}.  

Let $f,g\in{\mathcal C}^{\infty}({\mathfrak h}^*,K)$.
\begin{equation}   \label{deformed_rack_product}
(f\rhd_{\hbar}g)(\alpha)\,:=\,\sum_{r=0}^{\infty}\frac{\hbar^r}{r!}
\sum_{i_1,\ldots,i_r=1}^n\frac{\partial^r f}{\partial\alpha_{i_1}\dots
\partial\alpha_{i_r}}(0)\Big((\widetilde{\rm ad}_{i_1}\circ\ldots\circ
\widetilde{\rm ad}_{i_r})(g)\Big)(\alpha).
\end{equation} 
  
\begin{theo}  \label{compatibility_rack_structures}
For all $x,y\in{\mathfrak h}$, we have 
$$
e^{\hat{x}}\rhd_{\hbar} e^{\hat{y}}\,=
\,e^{\widehat{x\blacktriangleright_{\hbar} y}},
$$
where $\blacktriangleright_{\hbar}:{\mathfrak h}\times{\mathfrak h}
\to{\mathfrak h}$ is the Lie rack structure 
(\ref{EqDefStandardLieRackStructureOnLeibniz})
defined by exponentiating the adjoint action of the Leibniz algebra:
$$
x\blacktriangleright_{\hbar} y\,=\,e^{\hbar\,\ad_x}(y).
$$
\end{theo}
\begin{prf}
The proof of the theorem relies on the rack bialgebra structure 
on $\mathsf{S}({\mathfrak h})$ given by 
Definition \ref{definition_UARinfini} and is performed in the 
following lemmas:
\begin{lem}
The map ``hat'' ~$\hat{}:{\mathfrak h}\to{\mathcal C}^{\infty}({\mathfrak h}^*
,K)$ which sends $x\in{\mathfrak h}$ to the linear function $\hat{x}$ extends
to an injective morphism of commutative associative unital algebras
$\Psi:\mathsf{S}({\mathfrak h})\to{\mathcal C}^{\infty}({\mathfrak h}^*,K)$
such that 
$$
\Psi(x_1\bullet\ldots\bullet x_k)\,=\,\hat{x_1}\dots\hat{x_k}
$$
for all integers $k$ and all $x_1,\ldots,x_k\in{\mathfrak h}$. 
\end{lem}
\begin{prf}
This follows immediately from the freeness property of 
the algebra $\mathsf{S}({\mathfrak h})$.
\end{prf}  
\begin{lem}
The morphism $\Psi$ intertwines the adjoint actions 
$\widetilde{\rm ad}_{i}$ and $\ad_{e_i}^s$ (see eqn 
(\ref{adjoint_action_s})), i.e. for all $i=1,\ldots,n$, we have
$$
\widetilde{\rm ad}_{i}(\Psi(a))\,=\,\Psi(\ad_{e_i}^s(a))
$$
for all $a\in\mathsf{S}({\mathfrak h})$.
\end{lem}
\begin{prf}
Indeed, it is enough to show this for 
$x\in{\mathfrak h}\subset\mathsf{S}({\mathfrak h})$ as both adjoint actions 
are derivations. Now we have for $\alpha\in{\mathfrak h}^*$: 
\begin{eqnarray*}
\Psi(\ad_{e_i}^s(x))(\alpha) &=& \widehat{[e_i,x]}(\alpha) 
=  \alpha([e_i,x]) 
= \sum_{j,k=1}^n\,\alpha_k\,e^k([e_i,e_j])\,x_j\\
&=& \sum_{j,k=1}^n\,c_{ij}^k\,\alpha_k\,
\frac{\partial \hat{x}}{\partial\alpha_j}(\alpha) 
= \widetilde{\rm ad}_{i}(\Psi(a))(\alpha)
\end{eqnarray*}
\end{prf}
\begin{lem}
For all $b\in\mathsf{S}({\mathfrak h})$ and all 
$x_1,\ldots,x_r\in{\mathfrak h}$, we have
$$
\Psi(x_1\bullet\ldots\bullet x_r)\rhd_{\hbar}\Psi(b)\,=\,
\hbar^r\,\Psi((x_1\bullet\ldots\bullet x_r)\rhd b),
$$
where the left-hand $\rhd$ is the rack multiplication in the rack bialgebra 
$\mathsf{S}({\mathfrak h})$.
\end{lem}
\begin{prf}
First of all, note that by linearity it is enough to show this for 
$x_1,\ldots,x_r=e_{i_1},\ldots,e_{i_r}$ with $i_1,\ldots,i_r\in\{1,\ldots,n\}$. 
By eqn (\ref{adjoint_action_s}), we have 
$$
   (e_{i_1}\bullet\ldots\bullet e_{i_r})\rhd b\,=\,
\frac{1}{k!}\sum_{\sigma\in S_k}(\ad^s_{i_{\sigma(1)}}\circ\ldots\circ
\ad^s_{i_{\sigma(r)}})(b).
$$
Applying $\Psi$ gives then
\begin{eqnarray*}
\Psi((e_{i_1}\bullet\ldots\bullet e_{i_r})\rhd b)
    &=&
\frac{1}{k!}\sum_{\sigma\in S_k}\Psi((\ad^s_{i_{\sigma(1)}}\circ\ldots\circ
\ad^s_{i_{\sigma(r)}})(b)) \\
&=&
     \frac{1}{k!}\sum_{\sigma\in S_k}\widetilde{\rm ad}_{i_{\sigma(1)}}\circ
\ldots\circ\widetilde{\rm ad}_{i_{\sigma(r)}}(\Psi(b)),
\end{eqnarray*}
by the previous lemma. Now compute
$$
\frac{\partial^k\Big(\Psi(e_{i_1}\bullet\ldots\bullet e_{i_r})\Big)}
{\partial\alpha_{j_1}\ldots\partial\alpha_{j_k}}(0).
$$
This expression is non zero only if $k=r$ and 
$\{i_1,\ldots,i_k\}=\{j_1,\ldots,j_k\}$. In this case, the result is $1$.
One deduces the asserted formula.  
\end{prf}
Now we come back to the proof of the theorem. The assertion of the theorem
is the equality:
$$e^{\hat{x}}\rhd_{\hbar} e^{\hat{y}}\,=\,e^{\widehat{x\blacktriangleright_{\hbar} y}}.$$
Summing up the assertion of the previous lemma (taking $x_1=\ldots=x_r=x$), 
we obtain:
$$\sum_{r=0}^{\infty}\frac{1}{r!}\Psi(
\underbrace{x\bullet\ldots\bullet x}_{r\,\,\,{\rm times}})
\rhd_{\hbar}\Psi(b)\,=\,\Psi\left(\Big(\sum_{r=0}^{\infty}\frac{\hbar^r}{r!}
(\underbrace{x\bullet\ldots\bullet x}_{r\,\,\,{\rm times}})\Big)\rhd b\right),$$
and thus (as the rack product in $\mathsf{S}({\mathfrak h})$ is given 
by the adjoint action, using also that $\Psi$ is multiplicative)
\begin{equation}    \label{compatibility_rack_structures1}
e^{\hat{x}}\rhd_{\hbar}\Psi(b)\,=\,\Psi(e^{\hbar\,\ad_x}(b)).
\end{equation}   
This extends then to the asserted formula using that $e^{\hbar\,\ad_x}$ is 
an automorphism of $\mathsf{S}({\mathfrak h})$ (because it is the exponential 
of a derivation). 
\end{prf}

\begin{cor}  \label{main_result}
The above defined star-product induces the structure of a rack 
with respect to the product
$\rhd_{\hbar}$ on the set of exponential functions $U_{\mathfrak{h}}=\{E_{X}:=e^{\frac{i}{\hbar}X}\,|\, X\in{\mathfrak{h}}\}$ on ${\mathfrak h}^*$, and 
this star-product is equal (up to a sign) to the star-product in 
Theorem \ref{theorem_geometric_defoquant}, see also \cite{DW2013}.
\end{cor}

\begin{prf}
Via the formula of the theorem, the self-distributivity property of 
the rack product $\blacktriangleright$ in the rack bialgebra 
$\mathsf{S}({\mathfrak h})$ translates into the self-distributivity property 
of $\rhd_{\hbar}$ on the set of exponential functions.  
Since the star-product defined in \cite{DW2013} is a series of bidifferential
operators, and since such a series is uniquely determined by its values on
exponential functions, the present star-product coincides with the one found
in \cite{DW2013} thanks to the statement of the preceding theorem.
\end{prf}

\begin{rem}
Observe that the proof of the above theorem contains also an isomorphism
of commutative associative unital algebras between $\mathsf{S}({\mathfrak h})$
and the image of $\Psi$, i.e. the polynomial algebra generated by the 
$\hat{x_i}$, $i=1,\ldots,n$. Up to factors $\hbar^r$, $\Psi$ is by 
construction an isomorphism of rack bialgebras. This mirrors the 
relation between the universal enveloping algebra $\mathsf{U}({\mathfrak g})$
and the deformation quantization of $\mathsf{S}({\mathfrak g}^*)$ for a 
Lie algebra ${\mathfrak g}$, i.e. between $\mathsf{U}({\mathfrak g})$ and 
the Gutt star product on ${\mathfrak g}^*$, see \cite{Gut83}, \cite{DW2013}. 
\end{rem} 

\begin{rem}
 We would like to thank the referee for asking the interesting question of comparing the above star-product formula with
 the theory of Loday-Pirashvili, see \cite{LodPir} and
 \cite{Mos}: recall that they work in the
 so-called \emph{linear category} $\mathcal{LM}$ whose objects
 are diagrams $M\stackrel{\phi}{\to}N$ of vector spaces and
 linear maps with the obvious commuting squares of linear maps
 as morphisms. It can be seen as the category of complexes
 concentrated in degree $1$ ($M$) and $0$ ($N$) with differential $\phi$. There is the obvious tensor product
 of complexes (truncated in degree $\geq 2$) equipping
 $\mathcal{LM}$ with the structure of a symmetric monoidal category (with $0\to K$ as unit object and a signless flip
 as symmetric braiding). In this category, associative algebra
 objects will be dialgebras, and Lie objects can be seen as
 Leibniz algebras $\mathfrak{h}\to \overline{{\mathfrak h}}=\mathfrak{g}$,
 see \cite{LodPir} and \cite{Mos} for more details. The symmetric
 algebra generated by $\mathfrak{h}\to \overline{{\mathfrak h}}$
 will be $\mathfrak{h}\otimes \mathsf{S}(\mathfrak{g})
 \to \mathsf{S}(\mathfrak{g})$, and the universal enveloping
 algebra of $\mathfrak{h}\to \overline{{\mathfrak h}}$ will be
 $\mathfrak{h}\otimes \mathsf{U}(\mathfrak{g})
 \to \mathsf{U}(\mathfrak{g})$ where the linear maps are
 given by $x\otimes g\mapsto p(x)\bullet g$ and
 $x\otimes u\mapsto p(x)u$, respectively.\\
  For a possible rack-like
 star-product formula motivated by this category it seems to us not unreasonable to base it on the `underlying
 augmented rack-bialgebras' since this worked already well
 for the above rack-star-product (\ref{deformed_rack_product}) on smooth functions on
 $\mathfrak{h}^*$. In a slightly more general fashion,
 let $(B,\Phi_B,U(\mathfrak{g}),\ell)$ be an augmented rack
 bialgebra for a given Lie algebra $\mathfrak{g}$. We shall
 later specialize
 $B$ to $\mathsf{S}(\mathfrak{h})$ or its rack sub-bialgebra
 $\mathsf{S}(\mathfrak{h})_{(1)}=K\mathbf{1}\oplus \mathfrak{h}$
 mentioned in Theorem \ref{TLeibnizToUARInfinityOfLeibniz},
 part 3. Now, $U(\mathfrak{g})$ clearly acts `diagonally'
 (i.e. using its comultiplication)
 on the $K$-module $B\otimes U(\mathfrak{g})$ via
 $\ell$ on the first factor $B$ and via the adjoint action
 on the second factor $\mathsf{U}(\mathfrak{g})$. Calling this action $\ell^\otimes$, and defining
 $\Phi:B\otimes \mathsf{U}(\mathfrak{g})\to 
 \mathsf{U}(\mathfrak{g})$ by $\Phi(b\otimes u)=\Phi_B(b)u$
 we easily see that 
 $(B\otimes \mathsf{U}(\mathfrak{g}),\Phi, \mathsf{U}(\mathfrak{g}), \ell^\otimes)$ is an augmented
 rack bialgebra whose rack multiplication $\triangleright'$
 reads
  \begin{equation}\label{EqCompRackMultForBOtimesUg}
        (b\otimes u)\triangleright'(c\otimes v) 
        = \sum_{(b),(u)}
           \ell_{ \Phi_B(b^{(1)})u^{(1)}}(c)\otimes
               \mathrm{ad}_{\Phi_B(b^{(2)})u^{(2)}}(v)
  \end{equation}
  where $b,c\in B$ and $u,v\in\mathsf{U}(\mathfrak{g})$.\\
  The interesting cases are $B=\mathsf{S}(\mathfrak{h})$
  and its subcoalgebra $B=\mathsf{S}(\mathfrak{h})_{(1)}=K\mathbf{1}\oplus \mathfrak{h}$ for a Leibniz algebra $\mathfrak{h}$. The latter is important since 
  $(K\mathbf{1}\oplus \mathfrak{h})\otimes  \mathsf{U}(\mathfrak{g})$ is the universal enveloping
  algebra of $\mathfrak{h}$ as the left adjoint functor to
  the functor associating to any \emph{bar-unital dialgebra}
  its `commutator'-Leibniz algebra, see e.g.~\cite[Thm.2.8]{ABRW}.
  Note that the classical universal enveloping algebra of
  $\mathfrak{h}$ in the
  sense of Loday-Pirashvili is just the submodule $\mathfrak{h}\otimes  \mathsf{U}(\mathfrak{g})$
  of $(K\mathbf{1}\oplus \mathfrak{h})\otimes  \mathsf{U}(\mathfrak{g})$, and this is easily seen to be
  a rack-subalgebra (NOT a subcoalgebra) of
  $\mathsf{S}(\mathfrak{h})\otimes \mathsf{U}(\mathfrak{g})$
  with respect to the multiplication $\triangleright'$ in the above eqn
  (\ref{EqCompRackMultForBOtimesUg}). The primitive part
  of $\mathsf{S}(\mathfrak{h})\otimes \mathsf{U}(\mathfrak{g})$
  is known to be Leibniz algebra isomorphic to the hemi-semidirect product $\mathfrak{h}\oplus
  \mathfrak{g}$, see \cite[Thm.2.7,4.]{ABRW}, hence
  the rack bialgebra $\mathsf{S}(\mathfrak{h})\otimes \mathsf{U}(\mathfrak{g})$ is isomorphic to the rack bialgebra
  $\mathsf{S}(\mathfrak{h}\oplus \mathfrak{g})$. In the finite-dimensional case over $K=\mathbb{R}$ or
  $K=\mathbb{C}$ we can apply
  formula (\ref{deformed_rack_product}) to this hemi-semidirect
  product situation thus giving a star-product formula
  on the function algebra 
  $\mathcal{C}^\infty(\mathfrak{h}^*\times \mathfrak{g}^*,K)$
  reflecting the above multiplication (\ref{EqCompRackMultForBOtimesUg}),
  and clearly rack star-products on the subspaces of functions
  (at most) linear in $\mathfrak{h}^*$ which will then be
  rack star-product versions corresponding to the universal
  enveloping algebra 
  $\mathfrak{h}\otimes \mathsf{U}(\mathfrak{g})$ in the sense
  of Loday-Pirashvili.\\
  Surprisingly, there is a morphism of rack algebras
  $\Gamma:\mathsf{S}(\mathfrak{h})\to B\otimes \mathsf{U}(\mathfrak{g})$
  for $B=K\mathbf{1}\oplus\mathfrak{h}$, cf the constructions in \cite{Riv},
  \cite{RivWag}:
  firstly, the following general statement is very easy to check: given any two augmented
  rack bialgebras $(B,\Phi_B,H,\ell)$ and $(C,\Phi_C,H,\ell')$
  over the same cocommutative Hopf algebra $H$, then any $K$-linear map $\Gamma:B\to C$ intertwining the $H$-actions
  which satisfies $\Phi_C\circ \Gamma=\Phi_B$ will give
  a morphism of rack-algebras, 
  i.e.~$\Gamma(b\triangleright b')=\Gamma(b)\triangleright \Gamma(b')$
  for all $b,b'\in B$. Note that $\Gamma$ does not have to be a morphism of coalgebras, and will in general 
  NOT be a morphism of 
  rack bialgebras. More concretely, for
  $\Gamma:B=\mathsf{S}(\mathfrak{h})\to  C=(K\mathbf{1}\oplus\mathfrak{h})\otimes
   \mathsf{U}(\mathfrak{g})$ we make the following ansatz:
   \begin{equation}\label{EqDefGammaAnsatz}
   \Gamma= \big(
   (\mathbf{1}\epsilon+\mathrm{pr})
   \otimes (F_*(e^{(1)})\circ \Phi)\big)\circ \Delta.
   \end{equation}
   Here the maps $\Delta,\epsilon,\mathbf{1}$ give the usual 
   augmented coalgebra
   structure of $\mathsf{S}(\mathfrak{h})$, $\mathrm{pr}:
   \mathsf{S}(\mathfrak{h})\to \mathfrak{h}$ is the canonical
   projection, and
   $\Phi:\mathsf{S}(\mathfrak{h})\to \mathsf{U}(\mathfrak{g})$
   is the above map $\omega\circ \mathsf{S}(p)$, see eqn
   (\ref{EqDefPBWCircSp}). Next, $*$ denotes the convolution 
   multiplication on 
   $\mathrm{Hom}_K(\mathsf{U}(\mathfrak{g}),
   \mathsf{U}(\mathfrak{g}))$, and $e^{(1)}:\mathsf{U}(\mathfrak{g})
   \to \mathfrak{g}\subset \mathsf{U}(\mathfrak{g})$ is the 
   \emph{eulerian idempotent of the cocommutative bialgebra
   	$\mathsf{U}(\mathfrak{g})$}, 
   i.e.~$e^{(1)}=
   \ln_*\big(\mathbf{1}_{\mathsf{U}(\mathfrak{g})}
   \epsilon_{\mathsf{U}(\mathfrak{g})}
   +(\mathrm{id}_{\mathsf{U}(\mathfrak{g})}-
   \mathbf{1}_{\mathsf{U}(\mathfrak{g})}
   \epsilon_{\mathsf{U}(\mathfrak{g})})\big)$,
   see e.g.~\cite[Ch.4, p.139-141]{Loday}. Finally,
   $F(s)$ is a formal series with rational coefficients,
   and $F_*(e^{(1)})\in \mathrm{Hom}_K(\mathsf{U}(\mathfrak{g}),
   \mathsf{U}(\mathfrak{g}))$ denotes the convolution series
   where $s$ is replaced by $e^{(1)}$. Observe that
   the bialgebra structures of $\mathsf{S}(\mathfrak{h})$
   and of $\mathsf{U}(\mathfrak{g})$ are
   $\mathsf{U}(\mathfrak{g})$-module maps (recall that it is
   the adjoint action on $\mathsf{U}(\mathfrak{g})$) whence
   all convolutions of these structures are
   $\mathsf{U}(\mathfrak{g})$-module maps, hence also
   $e^{(1)}$ and its convolution series. It follows that
   $\Gamma$ is a $\mathsf{U}(\mathfrak{g})$-module map.
   Moreover, since 
   $e^{(1)}\circ \omega=\mathrm{pr}_\mathfrak{g}$ we get
   $e^{(1)}\circ\Phi=p\circ \mathrm{pr}$, and since
   $e^{*e^{(1)}}=\mathrm{id}_{\mathsf{U}(\mathfrak{g})}$,
   the second condition $\Phi_C\circ \Gamma=\Phi_B$ gives
   --using the surjectivity of $\Phi$--
   \[
      \big(\mathbf{1}_{\mathsf{U}(\mathfrak{g})}
      \epsilon_{\mathsf{U}(\mathfrak{g})}+e^{(1)}\big)
      * F_*(e^{(1)}) = \mathrm{id}_{\mathsf{U}(\mathfrak{g})}
      =  e^{*e^{(1)}} 
   \]
   hence with $F(s)=\frac{e^s}{1+s}$ we get the morphism of
   rack-algebras. Moreover there is a morphism $\Psi:
   \mathsf{S}(\mathfrak{h})^+\to 
   \mathfrak{h}\otimes \mathsf{U}
   (\mathsf{g})$ given by a similar ansatz: 
   in the above eqn (\ref{EqDefGammaAnsatz}) we just replace
   the expression
   $(\mathbf{1}\epsilon+\mathrm{pr})$ by $\mathrm{pr}$ and the
   series $F$ by the series $G(s)=\frac{e^s-1}{s}$.
\end{rem}  

\subsection{General deformation theory for rack bialgebras}

In this section, $(R,\Delta,\epsilon,\mu,\mathbf{1})$ is a \textbf{cocommutative} rack-bialgebra 
over a general commutative ring $K$,
and we use the notation $r\rhd s$ to denote the rack product $\mu(r\otimes s)$ of two elements 
$r$ and $s$ of $R$. In this subsection, we will often drop the symbol $\Sigma$ in Sweedler's 
notation of (iterated) 
comultiplications, so that the $n$-iterated comultiplication of $r$ in $R$ reads
\[
r^{(1)}\otimes\cdots\otimes r^{(n)}:=(\Delta\otimes \text{\rm Id}^{\otimes n-1})\circ \cdots\circ\Delta(r)
\]

Let $K_\hbar=K[[\hbar]]$ denote the $K$-algebra of formal power series in the indeterminate 
$\hbar$ with coefficients in $K$.
If $V$ is a vector space over $K$, $V_\hbar$ stands for V$[[\hbar]]$.
Recall that if $W$ is a $K$-module, a $K_\hbar$-linear morphism from $V_\hbar$ to $W_\hbar$ 
is the same as a power series in $\hbar$ with coefficients in $\text{\rm Hom}_{K}(V,W)$ via 
the canonical map
\[
\text{\rm Hom}_{K_\hbar}(V_\hbar,W_\hbar)\cong \text{\rm Hom}_{K}(V,W)_\hbar.
\]
This identification will be used without extra mention in the following.

\begin{defi}
A \textbf{formal deformation of the rack product $\mu$} is a formal power series $\mu_\hbar$
\[
\mu_\hbar:=\sum_{n\geq 0} \hbar^n \mu_n
\]
in $\text{\rm Hom}_{K}(R\otimes R,R)_\hbar$, such that
\begin{enumerate}
\item $\mu_0=\mu$,
\item $(R_\hbar,\Delta,\epsilon,\mu_\hbar,\mathbf{1})$ is a rack bialgebra over $K_\hbar$.
\end{enumerate}
\end{defi}

\begin{exem}
For a Leibniz algebra $\mathfrak{h}$, we have introduced in Definition 
\ref{definition_UARinfini} a cocommutative augmented rack bialgebra
$\mathsf{UAR}^{\infty}(\mathfrak{h})$. Furthermore, the rack star product defined in 
eqn (\ref{deformed_rack_product}), restricted to 
$S(\mathfrak{h})_\hbar$, is a deformation of the trivial rack product of $S(\mathfrak{h})$ 
given for all $r,s\in S(\mathfrak{h})$ by
\[
r \rhd s := \epsilon(r)s
\] 
The self-distributivity relation is shown in a way very similar to the proof of Theorem 
\ref{compatibility_rack_structures}, see eqn (\ref{compatibility_rack_structures1}). 
\end{exem}

As in the classical setting of deformation theory of associatice products, we will relate 
our deformation theory of rack products to cohomology. For this, let us first examine an
introductory example:

\begin{exem}  \label{example_infinitesimal_deformations}
Let $(R,\rhd)$ be a rack bialgebra, and suppose there exists a deformation 
$\rhd_{\hbar}=\rhd+\hbar\omega$ of $\rhd$. The new rack product $\rhd_{\hbar}$ 
should satisfy the self-distributivity identity, i.e. for all $a,b,c\in R$
$$a\rhd_{\hbar}(b\rhd_{\hbar}c)\,=\, (a^{(1)}\rhd_{\hbar}b)\rhd_{\hbar}(a^{(2)}\rhd_{\hbar}c)$$
To the order $\hbar^0$, this is only the self-distributivity relation for $\rhd$. But to order 
$\hbar^1$ (neglecting order $\hbar^2$ and higher), we obtain:
$$\omega(a,b\rhd c)+a\rhd\omega(b,c)\,=\,\omega(a^{(1)}\rhd b,a^{(2)}\rhd c)+
\omega(a^{(1)},b)\rhd(a^{(2)}\rhd c)+(a^{(1)}\rhd b)\rhd\omega(a^{(2)},c).$$
It will turn out that this is the cocycle condition for $\omega$ in the deformation complex
which we are going to define. More precisely, we will have
\begin{enumerate}
\item $d_{2,0}\omega(a,b,c)\,=\,\omega(a,b\rhd c)$,
\item $d_{1,1}\omega(a,b,c)\,=\,a\rhd\omega(b,c)$,
\item $d_{1,0}\omega(a,b,c)\,=\,\omega(a^{(1)}\rhd b,a^{(2)}\rhd c)$,
\item $d_{3}^2\omega(a,b,c)\,=\,\omega(a^{(1)},b)\rhd(a^{(2)}\rhd c)$,
\item $d_{2,1}\omega(a,b,c)\,=\,(a^{(1)}\rhd b)\rhd\omega((a^{(2)},c)$. 
\end{enumerate}
This may perhaps help to understand the general definition of the operators $d_{i,\mu}^n$
for $i=1,\ldots,n$ and $\mu\in\{0,1\}$ further down. 

On the other hand, the requirement that $\rhd_{\hbar}$ should be a morphism of coalgebras
(with respect to the undeformed coproduct $\triangle$ of $R$) means
$$\triangle\circ\rhd_{\hbar}\,=\,(\rhd_{\hbar}\otimes\rhd_{\hbar})\circ\triangle^{[2]}.$$
This reads for $a,b\in R$ to the order $\hbar$ (neglecting higher powers of $\hbar$) as
$$\omega(a,b)^{(1)}\otimes\omega(a,b)^{(2)}\,=\,\omega(a^{(1)},b^{(1)})\otimes(a^{(2)}\rhd b^{(2)})+
(a^{(1)}\rhd b^{(1)})\otimes\omega(a^{(2)},b^{(2)}).$$
This is exactly the requirement that $\omega$ is a coderivation along $\rhd=\mu$, to be defined below.
\end{exem}

Recall that $R$ being a rack bialgebra means in particular that $\mu:R^{\otimes 2}\to R$ 
is a morphism of coassociative coalgebras. For all positive integer $n$, let 
$\mu^n:R^{\otimes n}\to R$ be the linear map defined inductively by setting
\begin{itemize}
\item $\mu^1:=\text{\rm Id}:R\to R$,
\item $\mu^2:=\mu:R^{\otimes 2}\to R$,
\item $\mu^n:=\mu\circ (\mu^1\otimes \mu^{n-1})$, $n\geq 3$,
\end{itemize}
so that
\[
\mu^n(r_1,\cdots,r_n)\,=\,r_1\rhd (r_2\rhd (\cdots \rhd (r_{n-1}\rhd r_n)\cdots))
\] 
for all $r_1,\ldots,r_n$ in $R$. 

\begin{prop}\label{proposition-mu-n-coalg-morphism-formula}
  For all $n\geq 1$, the map $\mu^n$ is a morphism of coalgebras satisfying
  \begin{eqnarray}
\mu^{i}(r_1^{(1)},\!\cdots\!,r_{i-1}^{(1)},r_i)\rhd \mu^{n-1}(r_1^{(2)},\!\cdots\!,
r_{i-1}^{(2)},r_{i+1},\!\cdots\!,r_{n}) &=& \mu^{n}(r_1,\!\cdots\!,r_{n}), \label{eq:property-mu-n-rhd} \\
\mu^n(r_1,\cdots,r_{i-1},r_i^{(1)}\rhd r_{i+1},\cdots,r_i^{(n+1-i)}\rhd r_{n+1})
&=&\mu^{n+1}(r_1,\cdots,r_{n+1}) \nonumber   \\
&&   \label{eq:property-mu-n-rhd2}
\end{eqnarray}
for all positive integers $i$ and $n$ such that $1\leq i < n$ and for all $r_1$, ..., $r_{n}$ in $R$.
\end{prop}

\begin{prf}

\begin{itemize}
\item eqn (\ref{eq:property-mu-n-rhd}):
Let us show that the assertion of eqn (\ref{eq:property-mu-n-rhd}) is true for all $n$ and $i$
with $1\leq i < n$ by induction over $i$. Suppose
that the induction hypothesis is true and compute
\begin{eqnarray*}
\lefteqn{\mu^{i}(r_1^{(1)},\!\cdots\!,r_{i}^{(1)},r_{i+1})\rhd \mu^{n-1}(r_1^{(2)},\!\cdots\!,
r_{i}^{(2)},r_{i+2},\!\cdots\!,r_{n})} \\
&&\left(r_1^{(1)}\rhd\mu^{i}(r_2^{(1)},\!\cdots\!,r_{i}^{(1)},r_{i+1})\right)
\rhd\left(r_1^{(2)}\rhd\mu^{n-2}(r_2^{(2)},\!\cdots\!,r_{i}^{(2)},r_{i+2},\!\cdots\!,r_{n})\right),
\end{eqnarray*}
which gives, thanks to the self-distributivity relation in the rack algebra $R$,
\begin{eqnarray*}
\lefteqn{r_1\rhd\left(\mu^{i}(r_2^{(1)},\!\cdots\!,r_{i}^{(1)},r_{i+1})\rhd\mu^{n-2}(r_2^{(2)},\!\cdots\!,
r_{i}^{(2)},r_{i+2},\!\cdots\!,r_{n})\right)} \\
&&\,=\,r_1\rhd\mu^{n-1}(r_2,\!\cdots\!,r_{n})\,=\,\mu^{n}(r_1,\cdots,r_n),
\end{eqnarray*}
where we have used the induction hypothesis. This proves the assertion.

\item eqn (\ref{eq:property-mu-n-rhd2}): The assertion follows here again from an easy induction 
using the self-distributivity relation.
\end{itemize} 
\end{prf}

If $(C,\Delta_C)$ and $(D,\Delta_D)$ are two coassociative coalgebras and $\phi:C\to D$ is a 
morphism of coalgebras, we denote by $\text{\rm Coder}(C,V,\phi)$ the vector space of 
{\it coderivations from $C$ to $V$ along $\phi$}, i.e. the vector space of linear maps 
$f:C\to D$ such that
\[
\Delta_D\circ f\,=\,(f\otimes \phi+\phi\otimes f)\circ\Delta_C
\]

Let us note the following permanence property of coderivations along a map under partial 
convolution which will be useful in the proof of the following theorem. For
a coalgebra $A$, maps $f:A\otimes B\to V$
and $g:A\otimes C\to V$ and some product $\rhd:V\otimes V\to V$, the {\it partial convolution}
of $f$ and $g$ is the map $f\star_{\rm part}g:A\otimes B\otimes C\to V$ defined for all
$a\in A$, $b\in B$ and $c\in C$ by 
$$(f\star_{\rm part}g)(a\otimes b\otimes c)\,:=\,f(a^{(1)}\otimes b)\rhd g(a^{(2)}\otimes c).$$  

\begin{lem} \label{coderivation_lemma}
Let $A$, $B$, $C$ and $V$ be coalgebras, $V$ carrying a product $\rhd$ which is supposed to be a 
coalgebra morphism. 
Let $f:A\otimes B\to V$ be a coderivation along $\phi$ and $g:A\otimes C\to V$ be a coalgebra morphism.
Then the partial convolution $f\star_{\rm part}g$ is a coderivation along $\phi\star_{\rm part}g$.
\end{lem}
 
\begin{prf}
We compute for all $a\in A$, $b\in B$ and $c\in C$
\begin{eqnarray*}
\lefteqn{\triangle_V\circ(f\star_{\rm part}g)(a\otimes b\otimes c)=\triangle_V(f(a^{(1)}\otimes b)
\rhd g(a^{(2)}\otimes c))} \\
&& = (f(a^{(1)}\otimes b))^{(1)}\rhd(g(a^{(2)}\otimes c))^{(1)}\otimes
(f(a^{(1)}\otimes b))^{(2)}\rhd(g(a^{(2)}\otimes c))^{(2)} \\
&&= (f(a^{(1)}\otimes b))^{(1)}\rhd g(a^{(2)}\otimes c^{(1)})\otimes
(f(a^{(1)}\otimes b))^{(2)}\rhd g(a^{(3)}\otimes c^{(2)}) \\
&&=\phi(a^{(1)}\otimes b^{(1)})\rhd g(a^{(2)}\otimes c^{(1)})\otimes
f(a^{(3)}\otimes b^{(2)})\rhd g(a^{(4)}\otimes c^{(2)})  + \\
&&+ f(a^{(1)}\otimes b^{(1)})\rhd g(a^{(2)}\otimes c^{(1)})\otimes
\phi(a^{(3)}\otimes b^{(2)})\rhd g(a^{(4)}\otimes c^{(2)}) \\
&&=(\phi\star_{\rm part}g)(a^{(1)}\otimes b^{(1)}\otimes c^{(1)})\otimes 
(f\star_{\rm part}g)(a^{(2)}\otimes b^{(2)}\otimes c^{(2)}) + \\
&&+(f\star_{\rm part}g)(a^{(1)}\otimes b^{(1)}\otimes c^{(1)})\otimes
(\phi\star_{\rm part}g)(a^{(2)}\otimes b^{(2)}\otimes c^{(2)}) \\
&&= \big((\phi\star_{\rm part}g)\otimes(f\star_{\rm part}g)
+ (f\star_{\rm part}g)\otimes(\phi\star_{\rm part}g)\big)\circ\triangle_{A\otimes B\otimes C}(a\otimes b\otimes c).  
\end{eqnarray*}
\end{prf}

\begin{defi}
  The \textbf{deformation complex of $R$} is the graded vector space $C^*(R;R)$ defined in 
degree $n$ by
\[
C^n(R;R)\,:=\,\text{\rm Coder}(R^{\otimes n},R,\mu^n)
\]
endowed with the differential $d_R:C^*(R;R)\to C^{*+1}(R;R)$ defined in degree $n$ by
\[
d_R^n:=\sum_{i=1}^{n} (-1)^{i+1}(d_{i,1}^n-d_{i,0}^n)\:+\:(-1)^{n+1}d^n_{n+1}
\]
where the maps $d_{i,1}^n$ and $d_{i,0}^n$ are defined respectively by
\[
d_{i,1}^n\omega (r_1,\cdots,r_{n+1}):= \sum_{(r_1),\cdots,(r_i)}   
\mu^i(r_1^{(1)},\cdots,r_{i-1}^{(1)},r_i)\rhd \omega(r_1^{(2)},\cdots,r_{i-1}^{(2)},r_{i+1},\cdots,r_{n+1})
\]
and
\[
d_{i,0}^n\omega (r_1,\cdots,r_{n+1}):= \sum_{(r_i)}
\omega(r_1,\cdots,r_{i-1},r_i^{(1)}\rhd r_{i+1},\cdots,r_i^{(n+1-i)}\rhd r_{n+1})
\]
and $d_{n+1}^n$ by
\begin{eqnarray*}
\lefteqn{  d_{n+1}^n\omega (r_1,\cdots,r_{n+1})  }  \\ &:=& 
\sum_{(r_1),\cdots,(r_{n-1})}   
\omega(r_1^{(1)},\cdots,r_{n-1}^{(1)},r_{n})\rhd \mu^n(r_1^{(2)},\cdots,r_{n-1}^{(2)},r_{n+1})
\end{eqnarray*}
for all $\omega$ in $C^n(R;R)$ and $r_1,\ldots,r_{n+1}$ in $R$. 
\end{defi}

\begin{theo}  \label{rack_bialgebra_cohomology_complex}
$d_R$ is a well defined differential.
\end{theo}

\begin{prf}
That $d_R$ is well defined means that it sends coderivations to coderivations. 
It suffices to show that this is already true for all maps $d_{i,1}^n$, $d_{i,0}^n$ 
and $d_{n+1}^n$, which is the case. For this, we use Lemma \ref{coderivation_lemma}.
Indeed, a cochain $\omega\in C^n(R;R)$ is a coderivation along $\mu^n$. By Proposition
\ref{proposition-mu-n-coalg-morphism-formula}, $\mu^n$ is a coalgebra morphism. 
On the other hand, it is clear from the formula for $d_{i,1}^n$ that $d_{i,1}^n$ is a
partial convolution with respect to the first $i-1$ tensor labels of $\mu^i$ and $\omega$. 
Therefore the Lemma applies to give that the result is a coderivation along the partial 
convolution of $\mu^i$ and $\mu^n$, which is just $\mu^{n+1}$ again by Proposition 
\ref{proposition-mu-n-coalg-morphism-formula}. This shows that $d_{i,1}^n\omega$ belongs to 
$\text{\rm Coder}(R^{\otimes n},R,\mu^{n+1})$ as expected. 
The maps $d_{i,0}^n$ and $d_{n+1}^n$ can be treated in a similar way.

The fact that $d_R$ squares to zero is related to the so-called {\it cubical 
identities} satisfied by the maps $d_{i,1}$ and the maps $d_{i,0}$, namely
$$d^{n+1}_{j,\mu}\circ d^{n}_{i,\nu}\,=\, d^{n+1}_{i+1,\nu}\circ d^{n}_{j,\mu}\,\,\,\,\,
{\rm for}\,\,\,\,\,
j\leq i\,\,\,\,\,{\rm and}\,\,\,\,\,\mu,\nu\in\{0,1\},$$
and auxiliary identities which express the compatibility of the maps $d_{i,1}$ 
and $d_{i,0}$
with $d_{n+1}^n$, and an identity involving $d_{n+1}^n$ and $d_{n+2}^{n+1}$. 
One could call this kind of object an {\it augmented cubical vector space}.   

We will not show the usual cubical relations, i.e. those which do not refer to 
the auxiliary coboundary map $d_{n+1}^n$, because these are well-known to hold for rack cohomology,
see \cite{Cov}, Corollary 3.12, and our case is easily adapted from there. 
One possibility of adaptation (in case one works over the real or complex numbers) is to 
take a Lie rack, write its rack homology complex (with trivial
coefficients in the real or complex numbers), and to apply the functor
of point-distributions. 



Let us show that the two following extra relations involving the extra 
face $d_{n+1}^n$ hold:
\begin{equation}
  \label{eq:extra-relation-dnPlus1-di}
  d_{i,\mu}^{n+1}\circ d_{n+1}^n = d_{n+2}^{n+1}\circ d_{i,\mu}^n
\end{equation}
for all $1\leq i \leq n$ and $\mu$ in $\{0,1\}$ and 
\begin{equation}
  \label{eq:extra-relation-dnPlus2-dnPlus1}
  d_{n+1,0}^{n+1}\circ d_{n+1}^n = d_{n+2}^{n+1}\circ d_{n+1}^n + d_{n+1,1}^{n+1}\circ d_{n+1}^n
\end{equation}

Indeed, if $\omega$ is a $n$-cochain and $r_1$, ..., $r_{n+2}$ are elements in $R$, then
\begin{align*}
( d_{i,1}^{n+1}&\circ d_{n+1}^n\omega)(r_1,\!\cdots\!,r_{n+2})\\
=& \mu^i(r_1^{(1)},\!\cdots\!,r_{i-1}^{(1)},r_i)\rhd d_{n+1}^n\omega(r_1^{(2)},\!\cdots\!,
r_{i-1}^{(2)},r_{i+1},,\!\cdots\!,r_{n+2})\\
=&\mu^i(r_1^{(1)},\!\cdots\!,r_{i-1}^{(1)},r_i)\rhd \big(\omega(r_1^{(2)},\!\cdots\!,
r_{i-1}^{(2)},r_{i+1}^{(1)},,\!\cdots\!,r_n^{(1)},r_{n+1})\rhd\\
&\quad\rhd \mu^n(r_1^{(3)},\!\cdots\!,r_{i-1}^{(3)},r_{i+1}^{(2)},,\!\cdots\!,r_n^{(2)},
r_{n+2})\big)\\
\end{align*}
By Proposition \ref{proposition-mu-n-coalg-morphism-formula} and thanks to the 
self-distributivity of the rack product, this equality can be rewritten as
\begin{align*}
( d_{i,1}^{n+1}&\circ d_{n+1}^n\omega)(r_1,\!\cdots\!,r_{n+2})\\
=\,&\big(\mu^i(r_1^{(1)},\!\cdots\!,r_{i-1}^{(1)},r_i^{(1)})\rhd \omega(r_1^{(2)},\!\cdots\!,
r_{i-1}^{(2)},r_{i+1}^{(1)},,\!\cdots\!,r_n^{(1)},r_{n+1})\big)\rhd\\
&\quad\rhd \big( \mu^i(r_1^{(3)},\!\cdots\!,r_{i-1}^{(3)},r_i^{(2)})\rhd \mu^n(r_1^{(4)},
\!\cdots\!,r_{i-1}^{(4)},r_{i+1}^{(2)},,\!\cdots\!,r_n^{(2)},r_{n+2})\big)\\
=\,& d_{i,1}^{n}\omega(r_1^{(1)},\!\cdots\!,r_n^{(1)},r_{n+1})\rhd \mu^n(r_1^{(2)},\!\cdots\!,
r_n^{(2)},r_{n+2})\\
=\,& (d_{n+2}^{n+1}\circ d_{i,1}^n \omega)(r_1,\cdots,r_{n+2}),
\end{align*}
which proves that Relation (\ref{eq:extra-relation-dnPlus1-di}) holds when $\mu=1$. 
The case $\mu=0$ goes as follows:
\begin{align*}
  (d_{i,0}^{n+1}\circ &d_{n+1}^n \omega)(r_1,\cdots,r_{n+2})=\,d_{n+1}^n \omega(r_1,\cdots,
r_{i-1},r_i^{(1)}\rhd r_{i+1},\cdots,r_i^{(n+2-i)}\!\!\!\!\!\!\!\rhd r_{n+2})\\
=\,&\omega(r_1^{(1)},\!\cdots\!,r_{i-1}^{(1)},r_i^{(1)}\rhd r_{i+1}^{(1)},\!\cdots\!,r_i^{(n-i)}
\!\rhd r_n^{(1)},r_i^{(n+1-i)}\!\rhd r_{n+1})\:\rhd \\
&\quad\rhd \mu^n(r_1^{(2)},\!\cdots\!,r_{i-1}^{(2)},r_i^{(n+2-i)}\!\!\rhd r_{i+1}^{(2)},
\!\cdots\!,r_i^{(2n-2i+1)}\!\!\rhd r_n^{(2)},r_i^{(2n-2i+2)}\!\!\rhd r_{n+2})
\end{align*}
where we have used that the rack product is a morphism of coalgebras. 
Recall the following equation from
Proposition \ref{proposition-mu-n-coalg-morphism-formula}:
\[
\mu^n(s_1,\cdots,s_{i-1},s_i^{(1)}\rhd s_{i+1},\cdots,s_i^{(n+1-i)}\rhd s_{n+1})=\,\mu^{n+1}(s_1,
\cdots,s_{n+1})
\]
for all $s_1$, ..., $s_{n+1}$ in $R$ and $1\leq i\leq n$. This allows to rewrite the 
preceeding equality as
\begin{align*}
  (d_{i,0}^{n+1}\circ &d_{n+1}^n \omega)(r_1,\cdots,r_{n+2})\\
=\,&\omega(r_1^{(1)},\!\cdots\!,r_{i-1}^{(1)},r_i^{(1)}\rhd r_{i+1}^{(1)},\!\cdots\!,r_i^{(n-i)}
\!\rhd r_n^{(1)},r_i^{(n+1-i)}\!\rhd r_{n+1})\:\rhd \\
&\quad\rhd \mu^{n+1}(r_1^{(2)},\!\cdots\!,r_{i-1}^{(2)},r_i^{(n+2-i)}, r_{i+1}^{(2)},\!\cdots\!, 
r_n^{(2)}, r_{n+2})\\
=\,&d_{i,0}^n\: \omega(r_1^{(1)},\!\cdots\!,r_n^{(1)}, r_{n+1})\:\rhd\: \mu^{n+1}(r_1^{(2)},
\!\cdots\!, r_n^{(2)}, r_{n+2})\\
=\,&(d^{n+1}_{n+2}\circ d_{i,0}^{n})(r_1,\cdots,r_{n+2})
\end{align*}
which proves that (\ref{eq:extra-relation-dnPlus1-di}) holds when $\mu=0$. Relation 
(\ref{eq:extra-relation-dnPlus2-dnPlus1}) relies on the fact that cochains are 
coderivations. Indeed, 
\begin{align*}
 (d_{n+1,0}^{n+1}\circ& d_{n+1}^n \omega)(r_1,\cdots,r_{n+2}) =\,d_{n+1}^n\omega(r_1,\cdots,
r_n,r_{n+1}\rhd r_{n+2})\\
=\,&\omega(r_1^{(1)},\cdots,r_{n-1}^{(1)},r_n)\rhd\mu^n(r_1^{(2)},\cdots,r_{n-1}^{(2)},r_{n+1}
\rhd r_{n+2})\\
=\,&\omega(r_1^{(1)},\cdots,r_{n-1}^{(1)},r_n)\rhd\mu^{n+1}(r_1^{(2)},\cdots,r_{n-1}^{(2)},
r_{n+1}, r_{n+2})\\
=\,&\omega(r_1^{(1)},\cdots,r_{n-1}^{(1)},r_n)\rhd\Big(\mu^{n}(r_1^{(2)},\cdots,r_{n-1}^{(2)},
r_{n+1})\rhd \mu^{n}(r_1^{(3)},\cdots,r_{n-1}^{(3)}, r_{n+2})\Big)
\end{align*}
where we have used Proposition \ref{proposition-mu-n-coalg-morphism-formula} in the 
last equality. By self-distributivity of $\rhd$ and because $\omega$ is a coderivation, 
this gives
\begin{align*}
 (d_{n+1,0}^{n+1}\!\!\circ d_{n+1}^n \omega)(r_1,\!\!\cdots\!,r_{n+2})& =\,\big(\omega(
r_1^{(1)},\!\!\cdots\!,r_{n-1}^{(1)},r_n)^{(1)}\rhd\mu^{n}(r_1^{(2)},\!\!\cdots\!,
r_{n-1}^{(2)},r_{n+1})\big)\rhd\\
&\qquad \big(\omega(r_1^{(1)},\!\!\cdots\!,r_{n-1}^{(1)},r_n)^{(2)}\rhd \mu^{n}(r_1^{(3)},
\!\!\cdots\!,r_{n-1}^{(3)}, r_{n+2})\big)\\
=\,&\big(\omega(r_1^{(1)},\!\!\cdots\!,r_{n-1}^{(1)},r_n^{(1)})\rhd\mu^{n}(r_1^{(2)},
\!\!\cdots\!,r_{n-1}^{(2)},r_{n+1})\big)\rhd\\
&\quad \big(\mu^n(r_1^{(3)},\!\!\cdots\!,r_{n-1}^{(3)},r_n^{(2)})\rhd \mu^{n}(r_1^{(4)},
\!\!\cdots\!,r_{n-1}^{(4)}, r_{n+2})\big)\\
&+\big(\mu^n(r_1^{(1)},\!\!\cdots\!,r_{n-1}^{(1)},r_n^{(1)})\rhd\mu^{n}(r_1^{(2)},
\!\!\cdots\!,r_{n-1}^{(2)},r_{n+1})\big)\rhd\\
&\quad \big(\omega(r_1^{(3)},\!\!\cdots\!,r_{n-1}^{(3)},r_n^{(2)})\rhd \mu^{n}(r_1^{(4)},
\!\!\cdots\!,r_{n-1}^{(4)}, r_{n+2})\big)
\end{align*}
Applying Proposition \ref{proposition-mu-n-coalg-morphism-formula} again enables us to rewrite 
this last equality as
\begin{align*}
 (d_{n+1,0}^{n+1}\!\!\circ  d_{n+1}^n& \omega)(r_1,\!\!\cdots\!,r_{n+2})\\
=\,&\big(\omega(r_1^{(1)},\!\!\cdots\!,r_n^{(1)})\rhd\mu^{n}(r_1^{(2)},\!\!\cdots\!,
r_{n-1}^{(2)},r_{n+1})\big)\rhd  \mu^{n+1}(r_1^{(3)},\!\!\cdots\!,r_{n-1}^{(3)},r_n^{(2)},r_{n+2})\\
&+\mu^{n+1}(r_1^{(1)},\!\!\cdots\!,r_n^{(1)},r_{n+1})\rhd \big(\omega(r_1^{(2)},\!\!\cdots\!,
r_n^{(2)})\rhd \mu^{n}(r_1^{(3)},\!\!\cdots\!,r_{n-1}^{(3)}, r_{n+2})\big)\\
=\,&d_{n+1}^n\omega(r_1^{(1)},\cdots,r_n^{(1)},r_{n+1})\rhd  \mu^{n+1}(r_1^{(2)},\!\!\cdots\!,
r_n^{(2)},r_{n+2})\\
&+\mu^{n+1}(r_1^{(1)},\!\!\cdots\!,r_{n}^{(1)},r_{n+1})\rhd d_{n+1}^n\omega(r_1^{(2)},
\!\!\cdots\!,r_n^{(2)},r_{n+2})\\
=\,&\big((d_{n+2}^{n+1}\circ d_{n+1}^n+d_{n+1,1}^{n+1}\circ d_{n+1}^n)\: \omega \big)(r_1,\cdots,
r_{n+2})
\end{align*}
which proves (\ref{eq:extra-relation-dnPlus2-dnPlus1}). 

Let us show now how $d_R\circ d_R=0$ can be deduced from 
(\ref{eq:extra-relation-dnPlus1-di}), (\ref{eq:extra-relation-dnPlus2-dnPlus1}) 
and from the cubical relations. In degree $n$, we have
\begin{align*}
  d_R\circ d_R=\,&\big(\sum_{i=1}^{n+1}(-1)^{i+1} (d_{i,1}^{n+1}\!\!-d_{i,0}^{n+1})+
(-1)^{n+2}d_{n+2}^{n+1}\big) \circ\big(\sum_{i=1}^n(-1)^{i+1} (d_{i,1}^n\!\!-d_{i,0}^n)\\
&\qquad\qquad\qquad\qquad\qquad\qquad\qquad\qquad\qquad\qquad+(-1)^{n+1}d_{n+1}^{n}\big)\\
=\,&\sum_{i=1}^{n+1}\sum_{j=1}^n (-1)^{i+j} (d_{i,1}^{n+1}\circ d_{j,1}^n-d_{i,1}^{n+1}\circ 
d_{j,0}^n-d_{i,0}^{n+1}\circ d_{j,1}^n+d_{i,0}^{n+1}\circ d_{j,0}^n)\\
&+\sum_{i=1}^n(-1)^{n+i+1} (d_{n+2}^{n+1}\circ d_{i,1}^n-d_{n+2}^{n+1}\circ d_{i,0}^n-
d_{i,1}^{n+1}\circ d_{n+1}^{n}+d_{i,0}^{n+1} \circ d_{n+1}^{n})\\
&-d_{n+2}^{n+1}\circ d_{n+1}^n - d_{n+1,1}^{n+1}\circ d_{n+1}^n +d_{n+1,0}^{n+1}\circ d_{n+1}^n
\end{align*}
The first double sum is equal to zero thanks to the cubical relations, the second 
sum is zero thanks to relation (\ref{eq:extra-relation-dnPlus1-di}). Relation 
(\ref{eq:extra-relation-dnPlus2-dnPlus1}) implies that the last one vanishes. 
This shows that $d_R$ is indeed a differential and concludes the proof of the proposition.
\end{prf}

\begin{defi}
The cohomology of the deformation complex $(C^*(R;R),d_R)$ is 
called the adjoint cohomology of the rack bialgebra $R$ and is denoted by $H^*(R;R)$. 

\end{defi}

\begin{defi}
  An \textbf{infinitesimal deformation} of the rack product is a deformation of the rack 
product over the $K$-algebra of dual numbers $\bar{K}_\hbar:=\,K_\hbar/\raisebox{-.6ex}{$(\hbar^2)$}$,
 i.e. a linear map $\mu_1:R^{\otimes 2}\to R$ such that $\bar{R}_\hbar:=\, R\otimes \bar{K}_\hbar$ 
is a rack bialgebra over $\bar{K}_\hbar$ when equipped with $\mu_0+\hbar \mu_1$.\\
 Two infinitesimal deformations $\mu_0+\hbar\mu_1$ and $\mu_0+\hbar\mu'_1$ are said to 
be \textbf{equivalent} if there exists an automorphism $\phi:\bar{R}_\hbar\to \bar{R}_\hbar$ 
of the coalgebra of $(\bar{R}_\hbar,\Delta,\epsilon)$ of the form $\phi:=\,\id_R+\hbar\alpha$ 
such that 
\[
\phi\circ (\mu_0+\hbar \mu_1) =\,(\mu_0+\hbar \mu'_1)\circ\phi.
\]
\end{defi}

As usual, being equivalent is an equivalence relation and one has the following 
cohomological interpretation of the set of equivalence classes of infinitesimal 
deformations, denoted $Def(\mu_0,\bar{K}_\hbar)$:
\begin{prop}
 \[
Def(\mu_0,\bar{\mathbb{K}}_\hbar)\,=\,H^2(R;R)
 \]
The identificaton is obtained by sending each equivalence class 
$[\mu_0+\hbar \mu_1]$ in $Def(\mu_0,\bar{K}_\hbar)$ to the cohomology class 
$[\mu_1]$ in $H^2(R;R)$. 
\end{prop}

\begin{prf}
  One checks easily that the correspondence is well defined (if $\mu_0+\hbar \mu_1$ 
is an infinitesimal deformation, then $\mu_1$ is a $2$-cocycle, see Example
 \ref{example_infinitesimal_deformations}) and that it is 
bijective when restricted to equivalence classes.
\end{prf}

\begin{rem}
\begin{enumerate}
\item[(a)] The choice of taking coderivations in the deformation complex is explained as
follows: The rack product $\mu$ is a morphism of coalgebras, and we want to deform it 
as a morphism of coalgebras with respect to the fixed coalgebra structure we started with. 
Tangent vectors to $\mu$ in $\Hom_{\rm coalg}(C\otimes C,C)$ are exactly coderivations
along $\mu$.
This is the first step: Deformations as morphisms of coalgebras. Then as a second step,
we look for 1-cocycles, meaning that we determine those morphisms of coalgebras which give 
rise to rack bialgebra structures. The deformation complex in \cite{CCES08} takes 
into account also the possibility of deforming the coalgebra structure, and we recover 
our complex by restriction.  
\item[(b)] Given a Leibniz algebra ${\mathfrak h}$, there is a natural restriction map
from the cohomology complex with adjoint coefficients of ${\mathfrak h}$ to the
deformation complex of its augmented enveloping rack bialgebra $\mathsf{UAR}({\mathfrak h})$. 
The induced map in cohomology is not necessarily an isomorphism, as the abelian case shows.
Observe that the deformation complex of the rack bialgebra $K[R]$ for a rack $R$ does not 
contain the complex of rack cohomology for two reasons: First, this latter complex is
ill-defined for adjoint coefficients, and second, there are not enough coderivations as all 
elements are set-like. A way out for this last problem would be to pass to completions. 
\end{enumerate}
\end{rem}

\end{document}